\documentclass[a4paper,11pt]{article}
\usepackage{amsmath,amssymb,amsthm,graphics,latexsym,amsfonts,amsfonts}
\usepackage{fancyhdr}
\usepackage{algorithm, algpseudocode}
\usepackage[nooneline,flushleft]{caption2}
\usepackage{indentfirst}

\theoremstyle{plain}

\title{ \mbox{\Large Arc-connectivity and super arc-connectivity of  mixed Cayley digraph}}
\author{Yuhu Liu{\footnote{Corresponding author.\newline E-mail address:xjuliu@163.com(Y.H.Liu),  mjx@xju.edu.cn(J.Meng).  \newline The research is supported by NSFC ,FUKT and SPUX. }} \qquad Jixiang Meng \\
\small  College of Mathematics and System Sciences, Xinjiang
University \\ \small  Urumqi, Xinjiang, 830046, P.R.China }
\date{}

\begin{document}
\maketitle

\noindent{\bf Abstract} 

A digraph $X=(V, E)$ is max-~$\lambda$, if $\lambda(X)=\delta(X)$. A
digraph $X$ is super-$\lambda$ if every minimum  cut of X is either
the set of inarcs of some vertex or the set of outarcs of some
vertex. In this paper, we$^{'}$ll prove that for all but a few
exceptions, the strongly connected \emph{mixed Cayley digraphs} are
max-$\lambda$ and super$-\lambda$.


\section {\large Introduction}
Let $X=(V,E)$ be a digraph, where $V$ is a finite set and $E$ is an
irreflexive relation on $V$, thus $E$ is a set of ordered pairs
$(u,v)\in V\times V$ such that $u\neq v$, the elements of $V$ are
called the {\em vertices} or \emph{nodes} of $X$ and the elements of
$E$ are called the \emph{arcs} of $X$, arc $(u,v)$ is said to be an
\emph{inarc} of $v$ and an \emph{outarc} of $u$. If $u$ is a vertex
of $X$, then the \emph{outdegree} of $u$ in $X$ is the number
\emph{d}$_{X}^{+}$($u$) of arcs of $X$ originating at $u$ and the
\emph{indegree} of $u$  in $X$ is the number \emph{d}$_{X}^{-}$($u$)
of arcs of $X$ terminating at $u$. The minimum outdegree of $X$ is
$\delta^{+}$(X)=min{\{\emph{d}$_{X}^{+}$($u$) $|$ $u\in V$\}}, the
minimum indegree of $X$ is
$\delta^{-}$($X$)=min{\{\emph{d}$_{X}^{-}$($u$) $|$ $u\in V$\}}, we
denote by $\delta(X)$ the minimum of $\delta^{+}$($X$) and
$\delta^{-}$($X$).

An \emph{arc-disconnecting set} of $X$ is a subset $W$ of $E$ such
that $X$$\setminus$W=(V, E$\setminus$W) is not strongly connected.
An arc disconnecting set is \emph{minimal} if no proper subset of
$W$ is an arc disconnecting set of $X$ and is a \emph{minimum arc
disconnecting set} if no other arc disconnecting set has smaller
cardinality than $W$. The \emph{arc connectivity $\lambda$(X)} of a
nontrivial digraph $X$ is the cardinality of a minimum arc
disconnecting set of $X$.

The \emph{positive arc neighborhood} of a subset A of V is the set
$\omega_{X}^ {+}($A$)$ of all arcs which initiate at a vertex of A
and terminate at a vertex of $V$$\setminus$$A$. The \emph{negative
neighborhood} of subset $A$ of $V$ is the set $\omega_{X}^ {-}($A$)$
of all arcs which initiate in $V$$\setminus$$A$ and terminate in
$A$. Clearly $\omega_{X}^ {-}($A$)$=$\omega_{X}^ {+}($V$ \setminus
$A$)$. Arc neighborhoods of proper, nonempty subsets of $V$, often
called cuts, are clearly arc disconnecting sets.

An \emph{arc fragment} of $X$ is a proper, nonempty subset of $V$
whose positive or negative arc neighborhood has cardinality
$\lambda($X$)$.

We define a digraph $X$ to be \emph{super arc-connected}, or more
simply, \emph{super-$\lambda$}, if every minimum cut of $X$ is
either the set of inarcs of some vertex or the set of outarcs of
some vertex.

Let $X=(V, E)$ be a strongly connected digraph. An arc fragment of
least possible cardinality is called a \emph{$\lambda$-atom} of $X$
and a nontrivial arc fragment of least possible cardinality is
called a \emph{$\lambda$-superatom} of $X$. \\\\
\textbf{Definition 1.1} The \emph{reverse} digraph of digraph $X=(V,
E)$ is the digraph $X^{(r)}=(V,\{(v,u)$ $|$ $(u,v)\in E\})$, digraph
$X=(V,E)$ is \emph{symmetric} if $E$=$E^{(r)}$ and is
\emph{antisymmetric} if $E\bigcap E^{(r)}$=$\varnothing$.\\\\
\textbf{Definition 1.2} If $G$ is a group and $S$ is a subset of $G$
$\setminus{\{1_G\}}$, where $1_G$ is the identity of $G$. We define
the \emph{Cayley }\emph{digraph} \emph{Cay($G, S$)} to be the
digraph with vertices the elements of group G and arcs all pairs of
the form ($g$, $s\cdot$g) with $g$$\in$$G$ and $s$$\in$$S$. We
define a \emph{Cayley graph} to be a symmetric Cayley digraph. It
should be clear that a \emph{Cayley digraph Cay($G,S$)} is symmetric
if and only if the inverse of every element of $S$ is
again in $S$.\\\\
\textbf{Definition 1.3} Let $G$ be a group, $T_{0}$, $T_{1}\subseteq
 G$, the \emph{Bi-Cayley digraph} of $G$ with respect to $T_{0}$ and $T_{1}$
is defined as the bipartite digraph with vertex set $G\times\{0,
1\}$ and arc set \{(($g$, 0), ($t_{0}\cdot$$g$, 1)),
(($t_{1}\cdot$$g$, 1), ($g$, 0)) $|$ $g$$\in$$G$, $t_{0}\in$$T_{0}$,
$t_{1}\in$$T_{1}$ \}, denoted by $BD(G, T_{0}, T_{1})$.\\

J.X.Meng gives the definition of \emph{mixed Cayley digraph}. In
order to be convenient in this paper, we narrate it by another
way.\\\\
\textbf{Definition 1.4} Let $G$ be a finite group, $S_{0}$,
$S_{1}$$\subseteq$$G\setminus\{1_{G}\}$, $T_{0}$,
$T_{1}\subseteq$$G$. Define the \emph{mixed Cayley digraph}

$MD=G(G,S_{0},S_{1},T_{0},T_{1})$=Cay($G\times\{0\},S_{0}$)$\cup$Cay($G\times\{1\},S_{1}$)$\cup$$BD(G,T_{0},T_{1})$
as follows:

$1)$ $V(MD)=G\times\{0,1\}$, and let $X_{0}=G\times{\{0\}}$,
$X_{1}=G\times\{1\}$.

$2)$ $((g, i),(s_{i}\cdot g, i))\in E(MD)$, $g\in G,s_{i}\in S_{i}$,
for $i=0,1$.

$3)$ $((g,0),(t_{0}\cdot g, 1))\in E(MD)$, $((t_{1}\cdot$$g$,1),
(g,0))$\in E(MD)$ for $t_{0}\in$$T_{0}$, $t_{1}\in$$T_{1}$ and
$g\in$$G$.\\

So far, the research on the connectivity of the Cayley graph is
mainly focused on vertex connectivity, results on this subject are
referred to \cite{MJX,Tindell,Watkins}. The research on the
Bi-Cayley graph is primarily focused on its isomorphisms\cite{Lu},
few results, if any, are known on graphic properties of Bi-Cayley
graphs. The  results of Mixed Cayley graph are few. In \cite{Chen},
Chen and Meng point out that the Mixed Cayley graph also has high
connectivity. In this paper, we study the arc-connectivity of
strongly connected Mixed Cayley digraph, and we will prove that the
strongly connected Mixed Cayley digraphs are max-$\lambda$ and
super-$\lambda$ but a few exceptions.

We denote by $Aut(X)$ the automorphism group of $X$. The graph $X$
is said to be \emph{ vertex transitive} if $Aut(X)$ acts
transitively on $V(X)$, and to be \emph{edge transitive} if $Aut(X)$
acts transitively on $E(X)$. It is proved that these two kinds of
graphs usually have high connectivity. For instance, connected
vertex transitive graphs have maximum edge connectivity\cite{Mader},
and connected edge transitive graphs have maximum vertex
connectivity\cite{Tindell}.

For $a\in$$G$, the right multiplication $R^{'}(a)$: $g\rightarrow
ga$, $g\in G$, is clearly an automorphism of any Cayley digraph of
$G$. Let $R^{'}(G)$=\{$R^{'}(a)$: a$\in$G\}, then $R^{'}(G)$ is a
subgroup of \emph{the automorphism group} of any Cayley digraph. In
following proposition, we$^{'}$ll prove that $R(G)=\{R(a)|R(a):
(g,i)\rightarrow(ga,i),$ for $a, g\in G$ and i=0,1$\}$ is also a
subgroup of \emph{the automorphism group} of any \emph{mixed Cayley
digraph}.\\\\ \textbf{Proposition 1.5} Let
$X=MD(G,S_{0},S_{1},T_{0},T_{1})$, then

$(1)$ $R(G)\leqslant$\emph{Aut(X)}, thus \emph{Aut(X)} acts
transitively both on $X_{0}$ and $X_{1}$.

$(2)$ $d_{X}^{+}$((g,0))=$|T_{0}|+|S_{0}|$,
$d_{X}^{-}((g,0))=|T_{1}|+|S_{0}|$,

$d_{X}^{+}((g,1))=|T_{1}|+|S_{1}|$,
$d_{X}^{-}((g,1))=|T_{0}|+|S_{1}|$, for any $g\in G$.\\\\proof.
$(1)$ $((g_{1},i),(g_{2},i))\in E(X)\Leftrightarrow
g_{2}=s_{i}g_{1}$ for some $s_{i}\in
S_{i}$$\Leftrightarrow$$g_{2}a=s_{i}g_{1}a$$\Leftrightarrow$$((g_{1}a,i),\\
(g_{2}a,i))\in E(X)$$\Leftrightarrow$$R(a)((g_{1},i),(g_{2},i))\in
E(X)$ for $i=0,1$. \\$((g_{1},0),(g_{2},1))\in
E(X)$$\Leftrightarrow$$g_{2}=t_{0}g_{1}$ for some $t_{0}\in
T_{0}$$\Leftrightarrow$
$g_{2}a=t_{0}g_{1}a$$\Leftrightarrow$$((g_{1}a,0), (g_{2}a,1))\in
E(X)$$\Leftrightarrow$$R(a)((g_{1},0),(g_{2},1))\in
E(X)$.\\$((g_{2},1),(g_{1},0))\in
E(X)\Leftrightarrow$$g_{2}=t_{1}g_{1}$ for some $t_{1}\in
T_{1}$$\Leftrightarrow$$g_{2}a=t_{1}g_{1}a$$\Leftrightarrow$$((g_{2}a,1),(g_{1}a,0))\in
E(X)$$\Leftrightarrow$$R(a)((g_{2},1),(g_{1},0))\in E(X)$.\\So for
any $a\in G$, R(a) is an automorphism of the mixed Cayley digraph X,
thus $R(G)\leq \emph{Aut(X)}$, and since
$R(g_{1}^{-1}g_{2})((g_{1},i))=(g_{2},i)$ for any $g_{1},g_{2}\in
G$, \emph{Aut(X)} acts transitively both on $X_{0}$ and
$X_{1}$.\\$(2)$
$N^{+}((g,0))=\{\{T_{0}g\}\times\{1\}\}$$\cup$$\{\{S_{0}g\}\times\{0\}\}$,

$N^{-}((g,0))=\{\{T_{1}g\}\times\{1\}\}$$\cup$$\{\{S_{0}^{-1}g\}\times\{0\}\}$,

$N^{+}((g,1))=\{\{T_{1}^{-1}g\}\times\{0\}\}$$\cup$$\{\{S_{1}g\}\times\{1\}\}$

$N^{-}((g,1))=\{\{T_{0}^{-1}g\}\times\{0\}\}$$\cup$$\{\{S_{1}^{-1}g\}\times\{1\}\}$,\\
so we can get

$d_{X}^{+}((g,0))=|T_{0}|+|S_{0}|$,
$d_{X}^{-}((g,0))=|T_{1}|+|S_{0}^{-1}|=|T_{1}|+|S_{0}|$,

$d^{+}_{X}((g,1))=|T_{1}^{-1}|+|S_{1}|=|T_{1}|+|S_{1}|$
$d_{X}^{-}((g,1))=|T_{0}^{-1}|+|S_{1}^{-1}|=|T_{0}|+|S_{1}| $.
$\Box$

\section{Many results we need in this paper}\textbf{Proposition 2.1}\cite{Tindell} Let $X=(V,E)$ be a strongly
connected digraph and let $A$ and $B$ be positive(respectively,
negative) arc fragments of X such that $A\nsubseteq B$ and
$B\nsubseteq A$. If $A\cap B\neq\varnothing$ and $A\cup B\neq V$,
then each of the sets $A\cap B$, $A\cup B$, $A\setminus B$ and
$B\setminus A$ is an positive (respectively, negative) arc fragments of  $X$. $\Box$\\\\
\textbf{Corollary 2.2} Let $X=MD(G,S_{0},S_{1},T_{0},T_{1})$ be
strongly connected mixed Cayley digraph, If $\lambda(X)<\delta(X)$,
distinct positive(respectively, negative) $\lambda-$atoms are vertex
disjoint. $\Box$\\

An \emph{imprimitive block} for a group $\Phi$ of permutations of a
set $T$ is a proper, nontrivial subset $A$ of $T$ such that if
$\varphi\in \Phi$ then either $\varphi(A)=A$ or $\varphi(A)\cap A
=\varnothing$. \\\\ \textbf{Theorem 2.3}\cite{Tindell} Let $X=(V,E)$
be a graph or digraph and let $Y$ be the subgraph or subdigraph
induced by an imprimitive block $A$ of $X$. Then\\1. If $X$ is
vertex-transitive then so is $Y$.\\2. If $X$ is a strongly connected
arc-transitive digraph or a connected edge-transitive graph and $A$
is a proper subset of $V$, then $A$ is an independent subset of
$X$.\\3. If $X=Cay(G,S)$ and $A$ contains the identity of $G$, then
$A$ is a subgroup of $G$. $\Box$\\\\\textbf{Theorem
2.4}\cite{Tindell} If $X=(V,E)$ is a strongly connected digraph, but
not super$-\lambda$ and has $\delta(X)>2$, then distinct
positive(respectively, negative) $\lambda-$superatoms of $X$ are
vertex disjoint. $\Box$
\\\\\textbf{Theorem 2.5}\cite{Tindell} Every strongly connected vertex-transitive
digraph X satisfies $\lambda(X)=\delta(X)$. $\Box$
\section{Arc-connectivity of the mixed Cayley digraph }

 In this section, we$^{'}$ll prove that for all but a few exceptions, the mixed
Cayley digraph is max-$\lambda$. Clearly, if either $T_{0}$ or
$T_{1}$ is empty, $X=MD(G,S_{0},S_{1},T_{0},T_{1})$ isn$^{'}$t
strongly connected, so in following paper, we suppose that
$T_{0}\neq \varnothing$ and $T_{1}\neq \varnothing$.\\\\
\textbf{Proposition 3.1} Let $X=MD(G,S_{0},S_{1},T_{0},T_{1})$ be a
strongly connected mixed Cayley digraph and $A$ be a $\lambda$-atom.
If $\lambda(X)<\delta(X)$, then\\(1) $Y=X[A]$ is a strongly
connected subdigraph of $X$.\\(2) $|A|\geq \delta(X)+1$.\\\\Proof.
Without loss of generality, we suppose $A$ is a positive
$\lambda$-atom.\\(1) If $Y$ is not strongly connected, there exists
a proper subset $B$ of $A$ such that

$\omega_{Y}^{+}(B)=\varnothing$, so $\omega_{X}^{+}(B)\subseteq
\omega_{X}^{+}(A)$, \\thus

$|\omega_{X}^{+}(A)|=|\omega_{X}^{+}(B)|$ and $|B|<|A|$. \\It$^{'}$s
a contradiction. \\(2) Because
$\lambda(X)=|\omega_{X}^{+}(A)|\geq|A|(\delta(X)-|A|+1)$, if
$2\leq|A|\leq\delta(X)$, we can verify that

$|A|(\delta(X)-|A|+1)\geq\delta(X)$, \\thus when
$2\leq|A|\leq\delta(X)$, $\lambda(X)\geq\delta(X)$, it is a
contradiction.  $\Box$\\\\
\textbf{Lemma 3.2} Let $X=MD(G,S_{0},S_{1},T_{0},T_{1})$ be a
strongly connected mixed Cayley digraph and $A$ be a $\lambda$-atom.
If $\lambda(X)<\delta(X)$, then $|A\cap X_{i}|\geq2$  where
$X_{i}=\{$ $(g,i)$ $|$ $g\in G$ $\}$, for $i=0,1$.\\\\Proof. Without
loss of generality, suppose $A$ is a positive
$\lambda$-atom.\\\textbf{Claim 1} $A_{i}=A\cap X_{i}\neq\varnothing$
for $i=0,1$.\\If $A_{0}=\varnothing$ or $A_{1}=\varnothing$, then

$\lambda(X)=\omega_{X}^{+}(A)\geq
min\{|A||T_{0}|,|A||T_{1}|\}\geq|A|$,\\ thus by proposition 3.1,
$\lambda(X)\geq\delta(X)+1$, it is a contradiction.\\\textbf{Claim
2} $|A_{i}|=|A\cap X_{i}|\geq2$, for $i=0,1$.\\Suppose $|A_{0}|=1$,
then

$\lambda(X)=|\omega_{X}^{+}(A)|=\sum_{v\in A}d_{X}^{+}(v)-\sum_{v\in
A}d_{X[A]}^{+}(v)$=

$\sum_{v\in A_{0}}d_{X}^{+}(v)+\sum_{v\in
A_{1}}d_{X}^{+}(v)-\sum_{v\in A}d_{X[A]}^{+}(v)$. \\Because there
are at most $|T_{0}|+|T_{1}|$ arcs between $A_{0}$ and $A_{1}$,

$\lambda(X)\geq|T_{0}|+|S_{0}|+|T_{1}^{-1}||A_{1}|-(|T_{0}|+|T_{1}|)=|S_{0}|+(|A_{1}|-1)|T_{1}|$.\\
Since $X$ is strongly connected and $\lambda(X)<\delta(X)$, we have
$\delta(X)\geq2$. \\By proposition 3.1
$|A_{1}|=|A|-|A_{0}|=|A|-1\geq\delta(X)$, thus

$\lambda(X)\geq|S_{0}|+|T_{1}|=d_{X}^{-}((g,0))\geq\delta(X)$.\\It
is a contradiction. $\Box$\\\\\textbf{Lemma 3.3} Let
$X=MD(G,S_{0},S_{1},T_{0},T_{1})$ be a strongly connected mixed
Cayley digraph and $A=A_{0}\cup A_{1}$ be a $\lambda-$atom, where
$A_{i}=A\cap X_{i}=H_{i}\times\{i\}$, and $H_{i}\subseteq G$  for
$i=0,1$. Set $Y_{i}=X[A_{i}]$ be the subdigraph of $X$ induced by
$A_{i}$ for $i=0,1$. If $\lambda(X)<\delta(X)$, then\\$(1)$
\emph{Aut($Y_{i}$)} acts transitively on $A_{i}$ for $i=0,1$.\\(2)
If $A_{i}$ contains $(1_{G},i)$, then $H_{i}$ is a subgroup of $X$
for $i=0,1$.\\\\Proof. (1) By lemma 3.2, $A_{i}$ is nontrivial, for
any $(g_{1},i),(g_{2},i)\in A_{i}$, by proposition 1.5,
$R(g_{2}^{-1}g_{1})\in R(G)\leq\emph{Aut(X)}$. And it$^{'}$s easy to
verify that $R(g_{2}^{-1}g_{1})(A)$ is also a $\lambda-$atom, so
$R(g_{2}^{-1}g_{1})(A)=A$. Using proposition 1.5(1) and theorem 2.5,
we can deduce that $R(g_{2}^{-1}g_{1})(A_{i})=A_{i}$ for $i=0,1$. So
the restriction of $R(g_{2}^{-1}g_{1})$ on $A_{i}$ induces an
automorphism of $Y_{i}$, which maps $(g_{1},i)$ to $(g_{2},i)$.
Because $(g_{1},i)$ and $(g_{2},i)$ are two arbitrary vertices of
$A_{i}$, \emph{Aut($Y_{i})$} acts transitively on $A_{i}$ for
$i=0,1$.\\(2) By lemma 3.2, $|A_{i}|\geq2$. Then for any arbitrary
vertex $(g,i)\in A_{i},R(g^{-1})((g,i))=(1_{G},i)$, so
$R(g^{-1})(A)=A$, it means that $Ag^{-1}=A$, so $A_{i}g^{-1}=A_{i}$,
thus we get that $hg^{-1}\in H_{i}$, for any $h,g\in H_{i}$, so
$H_{i}$ is a subgroup of $G$ for $i=0,1$. $\Box$\\\\\textbf{Lemma
3.4} Let $X=MD(G,S_{0},S_{1},T_{0},T_{1})$ be a strongly connected
mixed Cayley digraph, and $A$ be a $\lambda-$atom of $X$, let
$A_{i}=A\cap X_{i}$, then if $\lambda(X)<\delta(X)$, we have
that\\(1) $V(X)$ is a disjoint union of distinct
positive(respectively,negative) $\lambda-$atoms of $X$.\\(2)
$|A_{0}|=|A_{1}|$.\\\\Proof. $(1)$ Without loss of generality, set
$A$ be a positive $\lambda-$atom, by proposition 1.5, \emph{Aut(X)}
acts transitively both on $X_{0}$ and $X_{1}$. Because
$\lambda(X)<\delta(X)$, from theorem 2.5, $X$ isn$^{'}$t vertex
transitive. Thus $X$ has exactly two orbits $X_{0}$ and $X_{1}$, by
lemma 3.2, $|A_{i}|\geq2$, so at least one vertex of $X_{i}$ lines
in $A$ respectively. So every vertex of $X$ lines in a positive
$\lambda-$atom. By corollary 2.2, $V(X)$ is a disjoint union of
distinct positive $\lambda-$atoms. \\(2) Let
$V(X)=\cup_{i=1}^{k}\varphi_{i}(A)$, where $\varphi_{i}\in
\emph{Aut(X)}$ such that $\varphi_{i}(A)\cap
\varphi_{j}(A)\neq\emptyset$ if and only if $i=j$, then
$X_{i}=\cup_{i=1}^{k}\varphi_{i}(A_{i})$. Since $|X_{0}|=|X_{1}|$,
we have $|A_{0}|=|A_{1}|$. $\Box$\\\\\textbf{Lemma 3.5} Let
$X=MD(G,S_{0},S_{1},T_{0},T_{1})$ be a strongly connected mixed
Cayley digraph with $\lambda(X)<\delta(X)$ and $A=A_{0}\cup A_{1}$
be a $\lambda-$atom, where $A_{i}=A\cap X_{i}=H_{i}\times\{i\}$ and
$H_{i}\subseteq G$ for $i=0,1$. Then\\$(1)$ If $(1_{G},0)\in A_{0}$,
then $H_{1}=t_{i}H_{0}$ for some $t_{i}\in T_{i}$, furthermore,

$X_{0}=\cup_{j=1}^{k}(H_{0}g_{j})\times\{0\}$,

$X_{1}=\cup_{j=1}^{k}(t_{i}H_{0}g_{j})\times\{1\}$, \\where
$R(g_{j})(A)\cap R(g_{l})(A)\neq\varnothing$ if and only if $j=l$
for $1\leq j,l \leq k$. \\(2) If $(1_{G},1)\in A_{1}$, then
$H_{0}=t_{i}^{-1}H_{1}$ for some $t_{i}\in T_{i}$, furthermore,

$X_{0}=\cup_{j=1}^{k}(t_{i}^{-1}H_{1}g_{j})\times\{0\}$,

$X_{1}=\cup_{j=1}^{k}(H_{1}g_{j})\times\{1\}$, \\where
$R(g_{j})(A)\cap R(g_{l})(A)\neq\varnothing$ if and only if $j=l$
for $1\leq j,l \leq k$.\\\\Proof. (1) Since $(1_{G},0)\in A_{0}$ and
$X[A]$ is strongly connected by proposition 3.1, there must exist at
least an element $t_{i}\in T_{i}$ such that $t_{i}\in H_{1}$. If
$(1_{G},0)\in A_{0}$, $H_{0}\leq G$. \\Then for any $h_{0}\in
H_{0}$, $R(h_{0})(A)=A$, since $H_{0}h_{0}=H_{0}$. \\Thus for any
$h_{0}\in H_{0}, H_{1}h_{0}=H_{1}$, so $H_{1}H_{0}=H_{1}$. \\And
because $t_{i}\in H_{1}$ and $|H_{0}|=|H_{1}|$, we have that
$H_{1}=t_{i}H_{0}$. \\Since $H_{0}\leq G$, we get that
$G=\cup_{j=1}^{k}(H_{0}g_{j})$, where $g_{1}=1_{G}$ and
$H_{0}g_{j}\cap H_{0}g_{l}\neq\varnothing$ if and only if $j=l$ for
$1\leq j,l\leq k$. Therefore,

 $V(X)=\cup_{j=1}^{k}R(g_{j})(A)$.\\So
$X_{0}=\cup_{j=1}^{k}R(g_{j})(A_{0})=\cup_{j=1}^{k}(H_{0}g_{j})\times\{0\}$,

$X_{1}=\cup_{j=1}^{k}R(g_{j})(A_{1})=\cup_{j=1}^{k}(H_{1}g_{j})\times\{1\}$=$\cup_{j=1}^{k}(t_{0}H_{0}g_{j})\times\{1\}$.\\(2)
It is similar to (1). $\Box$\\\\\textbf{proposition 3.6} Let
$X=MD(G,S_0,S_1,T_0,T_1)$ be a strongly connected mixed Cayley
digraph with $\lambda(X)<\delta(X)$, and let $A$ be a $\lambda-$atom
of $X$, and let $Y=X[A]$, then \emph{Aut(Y)} acts transitively both
on $A_{0}$ and $A_1$, where $A_i=A \cap X_i$ for $i=0,1$.\\\\Proof.
It is clearly true from proposition 1.5, corollary 2.2 and lemma
3.3.
$\Box$\\

Set $H=Y\setminus\{E(Y_0)\cup E(Y_1)\}$ where $Y_i=X[A_i],$ $i=0,
1$, then from lemma 3.3 and proposition 3.6, we can get

$d_{H}^{+}((g,0)),d_{H}^{-}((g,0)),d_{H}^{+}((g,1)),d_{H}^{-}((g,1))$,$d_{Y_i}^{+}(g,i)$
and $d_{Y_i}^{-}(g,i)$ are constant respectively. Furthermore,

$d_{H}^{+}((g,0))=d_{H}^{-}((g,1))$,
$d_{H}^{-}((g,0))=d_{H}^{+}((g,1))$ and
$d_{Y_i}^{+}(g,i)$=$d_{Y_i}^{-}(g,i)$. \\So we set
$d_{H}^{+}((g,0))=d_{H}^{-}((g,1))=p$,
$d_{H}^{-}((g,0))=d_{H}^{+}((g,1))=q$ and $Y_i$ is $r_i$ regular
digraph.

If $X$ is a strongly connected mixed Cayley digraph with
$\lambda(X)<\delta(X)$, from lemma 3.4, $V(X)$ is the union of
distinct positive (respectively, negative) $\lambda-$atoms of $X$.
Set $A$ is a $\lambda-$atom of $X$ and $A_i=A\cap X_i$, for i=0, 1.
Now we introduce a class of digraphs consisting of the following
eight classes of digraphs, denoted by $\Gamma$,

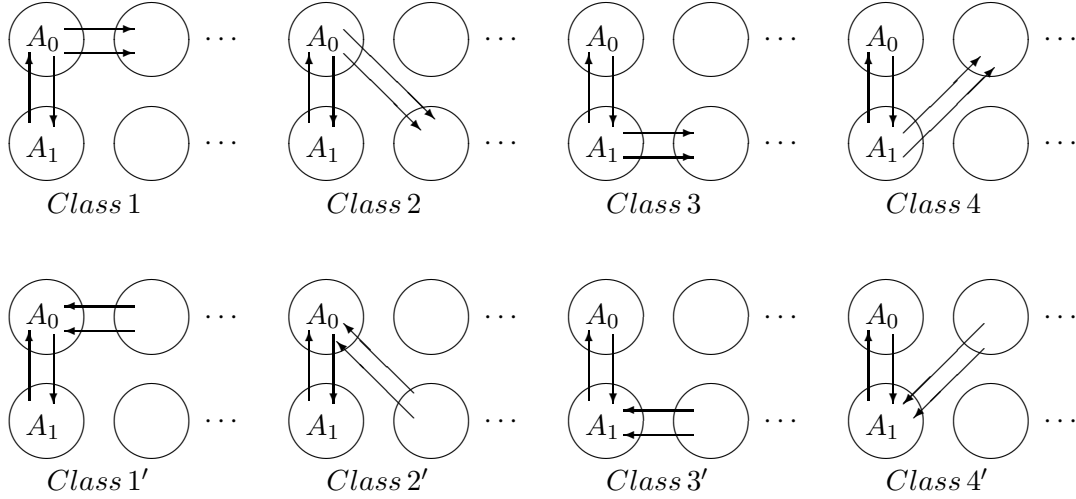
\begin{figure}[h,t]
\setlength{\unitlength}{0.46mm}
\begin{center}
\begin{picture}(290,136)

\put(0,10){\circle{20}}
\put(0,40){\circle{20}}
\put(0,90){\circle{20}}
\put(0,120){\circle{20}}

\put(30,10){\circle{20}}
\put(30,40){\circle{20}}
\put(30,90){\circle{20}}
\put(30,120){\circle{20}}

\put(80,10){\circle{20}}
\put(80,40){\circle{20}}
\put(80,90){\circle{20}}
\put(80,120){\circle{20}}

\put(110,10){\circle{20}}
\put(110,40){\circle{20}}
\put(110,90){\circle{20}}
\put(110,120){\circle{20}}

\put(160,10){\circle{20}}
\put(160,40){\circle{20}}
\put(160,90){\circle{20}}
\put(160,120){\circle{20}}

\put(190,10){\circle{20}}
\put(190,40){\circle{20}}
\put(190,90){\circle{20}}
\put(190,120){\circle{20}}

\put(240,10){\circle{20}}
\put(240,40){\circle{20}}
\put(240,90){\circle{20}}
\put(240,120){\circle{20}}

\put(270,10){\circle{20}}
\put(270,40){\circle{20}}
\put(270,90){\circle{20}}
\put(270,120){\circle{20}}

\put(0,-10){$Class\,1'$}
\put(0,70){$Class\,1$}

\put(80,-10){$Class\,2'$}
\put(80,70){$Class\,2$}

\put(160,-10){$Class\,3'$}
\put(160,70){$Class\,3$}

\put(240,-10){$Class\,4'$}
\put(240,70){$Class\,4$}

\put(-6,6){$A_1$}
\put(-6,38){$A_0$}
\put(45,38){$\cdots$}
\put(45,8){$\cdots$}

\put(-6,86){$A_1$}
\put(-6,118){$A_0$}
\put(45,118){$\cdots$}
\put(45,88){$\cdots$}

\put(74,6){$A_1$}
\put(74,38){$A_0$}
\put(125,38){$\cdots$}
\put(125,8){$\cdots$}

\put(74,86){$A_1$}
\put(74,118){$A_0$}
\put(125,118){$\cdots$}
\put(125,88){$\cdots$}

\put(154,6){$A_1$}
\put(154,38){$A_0$}
\put(205,38){$\cdots$}
\put(205,8){$\cdots$}

\put(154,86){$A_1$}
\put(154,118){$A_0$}
\put(205,118){$\cdots$}
\put(205,88){$\cdots$}

\put(234,6){$A_1$}
\put(234,38){$A_0$}
\put(285,38){$\cdots$}
\put(285,8){$\cdots$}

\put(234,86){$A_1$}
\put(234,118){$A_0$}
\put(285,118){$\cdots$}
\put(285,88){$\cdots$}

\put(-5,16){\vector(0,1){20}}
\put(2,35){\vector(0,-1){20}}
\put(-5,96){\vector(0,1){20}}
\put(2,115){\vector(0,-1){20}}

\put(75,16){\vector(0,1){20}}
\put(82,35){\vector(0,-1){20}}
\put(75,96){\vector(0,1){20}}
\put(82,115){\vector(0,-1){20}}

\put(155,16){\vector(0,1){20}}
\put(162,35){\vector(0,-1){20}}
\put(155,96){\vector(0,1){20}}
\put(162,115){\vector(0,-1){20}}

\put(235,16){\vector(0,1){20}}
\put(242,35){\vector(0,-1){20}}
\put(235,96){\vector(0,1){20}}
\put(242,115){\vector(0,-1){20}}

\put(25,36){\vector(-1,0){20}}
\put(25,43){\vector(-1,0){20}}
\put(5,116){\vector(1,0){20}}
\put(5,123){\vector(1,0){20}}

\put(185,6){\vector(-1,0){20}}
\put(185,13){\vector(-1,0){20}}
\put(165,86){\vector(1,0){20}}
\put(165,93){\vector(1,0){20}}

\put(85,116){\vector(1,-1){22}}
\put(85,123){\vector(1,-1){26}}

\put(245,86){\vector(1,1){26}}
\put(245,93){\vector(1,1){22}}

\put(105,11){\vector(-1,1){22}}
\put(105,18){\vector(-1,1){20}}
\put(268,31){\vector(-1,-1){20}}
\put(268,38){\vector(-1,-1){23}}
\end{picture}
\end{center}
\end{figure}
where $|A_0|=|A_1|<\delta(X)$ and Class 1 satisfies

$|S_0|-r_0=1,|T_0|-p=0,|S_1|-r_1=0$ and $|T_1|-q=0$.\\The Class 2
satisfies

$|S_0|-r_0=0,|T_0|-p=1,|S_1|-r_1=0$ and $|T_1|-q=0$.\\The Class 3
satisfies

$|S_0|-r_0=0,|T_0|-p=0,|S_1|-r_1=1$ and $|T_1|-q=0$. \\The Class 4
satisfies

$|S_0|-r_0=0,|T_0|-p=0,|S_1|-r_1=0$ and $|T_1|-q=1$.\\\\the Class
$1^{'}$ satisfies

$|S_0|-r_0=1,|T_0|-p=0,|S_1|-r_1=0$ and $|T_1|-q=0$.\\The Class
$2^{'}$ satisfies

$|S_0|-r_0=0,|T_0|-p=0,|S_1|-r_1=0$ and $|T_1|-q=1$.\\The Class
$3^{'}$ satisfies

$|S_0|-r_0=0,|T_0|-p=0,|S_1|-r_1=1$ and $|T_1|-q=0$. \\The Class
$4^{'}$ satisfies

$|S_0|-r_0=0,|T_0|-p=1,|S_1|-r_1=0$ and $|T_1|-q=0$.

Clearly, the  Class 1 and the  Class 3 are equivalent to the
Class $1^{'}$ and the  Class $3^{'}$ respectively. And we can also
easily prove that the  Class 2 and the  Class 4 are equivalent to the  Class $4^{'}$ and Class $2^{'}$ respectively.\\
\\\textbf{Theorem 3.7} Let $X=MD(G,S_0,S_1,T_0,T_1)$ be a strongly
connected mixed Cayley digraph. Then $X$ is not max$-\lambda$ if and
only if $X$ belongs to the class of digraphs $\Gamma$.\\\\Proof.
Necessity. If $\lambda(X)<\delta(X)$, by proposition 3.6, we set
that \\$d_{H}^{+}((g,0))=d_{H}^{-}((g,1))=p$ ,
$d_{H}^{-}((g,0))=d_{H}^{+}((g,1))=q$, and $Y_i$ is $r_i-$regular
digraph. Let $A$ be a $\lambda-$atom.\\\textbf{1}.
When $A$ is a positive $\lambda-$atom, then \\
$\lambda(X)=|\omega_X^+(A)|=|A_0|(|S_0|-r_0+|T_0|-p)+|A_1|(|S_1|-r_1+|T_1|-q)$.
\\Since $|A|\geq \delta(X)+1$ and $|A_0|+|A_1|\geq \delta(X)+1$, we
have $|A_0|=|A_1|> \delta(X)/ 2$. \\So
$\lambda(X)=|\omega_X^+(A)|<\delta(X)$ is true only if one of the
following conditions holds.\\\textbf{Case 1} $|S_0|-r_0+|T_0|-p=1$
and $|S_1|-r_1+|T_1|-q=0$.

\textbf{Subcase 1.1 }$|S_0|-r_0=1,|T_0|-p=0,|S_1|-r_1=0$ and
$|T_1|-q=0$,

clearly, under this subcase $X$ is Class 1.

\textbf{Subcase 1.2} $|S_0|-r_0=0,|T_0|-p=1,|S_1|-r_1=0$ and
$|T_1|-q=0$,

clearly, under this subcase $X$ is Class 2.
\\\textbf{Case 2} $|S_0|-r_0+|T_0|-p=0$ and $|S_1|-r_1+|T_1|-q=1$.

\textbf{Subcase 2.1} $|S_0|-r_0=0,|T_0|-p=0,|S_1|-r_1=1$ and
$|T_1|-q=0$,

clearly, under this subcase $X$ is Class 3.

\textbf{Subcase 2.2} $|S_0|-r_0=0,|T_0|-p=0,|S_1|-r_1=0$ and
$|T_1|-q=1$,

clearly, under this subcase $X$ is Class 4.
\\\textbf{2}. When $A$ is a
negative $\lambda-$atom, then\\
$\lambda$(X)=$|\omega_X^{-}(A)|=|A_0|(|S_0|-r_0+|T_1|-q)+|A_1|(|S_1|-r_1+|T_0|-p)$.\\
Since $|A|\geq \delta(X)+1$ and $|A_0|=|A_1|$, we have that
$|A_0|=|A_1|>\delta(X)/2$. \\So if
$\lambda$(X)=$|\omega_X^{-}(A)|<\delta(X)$, one of the following
conditions holds,\\\textbf{Case $1^{'} $} $|S_0|-r_0+|T_1|-q=1$ and
$|S_1|-r_1+|T_0|-p=0$.

\textbf{Subcase $1^{'}.1$} $|S_0|-r_0=1$,$|T_1|-q=0$,$|S_1|-r_1=0$
and $|T_0|-p=0$,

clearly, under this subcase $X$ is Class $1^{'}$.

\textbf{Subcase $1^{'}.2$} $|S_0|-r_0=0$,$|T_1|-q=1$,$|S_1|-r_1=0$
and $|T_0|-p=0$,

clearly, under this subcase $X$ is Class $2^{'}$.\\\textbf{Case
$2^{'} $} $|S_0|-r_0+|T_1|-q=0$ and $|S_1|-r_1+|T_0|-p=1$.

\textbf{Subcase $2^{'}.1$} $|S_0|-r_0=0$,$|T_1|-q=0$,$|S_1|-r_1=1$
and $|T_0|-p=0$,

clearly, under this subcase $X$ is Class $3^{'}$.

\textbf{Subcase $2^{'}.2$} $|S_0|-r_0=0$,$|T_1|-q=0$,$|S_1|-r_1=0$
and $|T_0|-p=1$,

clearly, under this subcase $X$ is Class $4^{'}$.\\Sufficency, it is clearly true. $\Box$\\
\\\textbf{Proposition 3.8} $X=MD(G,S_0,S_1,T_0,T_1)$ is a strongly connected mixed
Cayley digraph, $X$ is Class 1 or Class $1^{'}$ if and only if
\\(1) There exists a non-empty proper subgroup $H$ of $G$ and $S_0$
contains an element $s_0$ such that

$<S_0\cup \{1_G\}\setminus \{s_0\}>\leq H$ and $|H|<\delta(X)$, and
\\(2) There is an element $t_0\in T_0$ such that

 $G_1=<S_1>\leq
t_0Ht_0^{-1},T_1^{-1}t_0\subseteq H$ and $t_0^{-1}T_0\subseteq
H$.\\\\Proof. Necessity. Because  Class 1 is equivalent to Class
$1^{'}$, without loss of generality, we set $X$ is Class 1, then
Assume $(1_{G},0)\in A_0$, by lemma 3.3, $H_0\leq G$. Let $H=H_0$,
then under this situation we can achieve the following results
easily,\\(i)
$\lambda(X)=|\omega_X^{+}(A)|=|A_0|=|H_0|=|H|<\delta(X)$, since
$|S_0|-r_0=1$, $|T_0|-p=0$, $|S_1|-r_1=0$, and $|T_1|-q=0$,\\(ii)
$<S_0\cup\{1_G\}\setminus\{s_0\}>\leq H_0=H$, since $|S_0|-r_0=1$.\\
By proposition 3.5, $H_1=t_0H_0$ for some $t_0\in T_0$ and
$X_1=\cup_{i=1}^{k}(t_0H_0g_i)\times\{1\}$, where $t_0H_0g_i\cap
t_0H_0g_j\neq \varnothing$ if and only if $i=j$ for $1\leq i,j\leq
k$. Assume that $(1_G,1)\in (t_0H_0g_s)\times\{1\}$, then we can
deduce that $t_0H_0g_s\leq G$ and $g_s=h_0^{-1}t_0^{-1}$, where
$h_0\in H_0$.\\Since $|S_1|-r_1=0$, we get $G_1\leq
t_0H_0g_s=t_0H_0h_0^{-1}t_0^{-1}=t_0H_0t_0^{-1}=t_0Ht_0^{-1}$.\\Since
$|T_0|-p=0$ and $|T_1|-q=0$, we have that $T_0H_0\subseteq H_1$ and
$T_1^{-1}H_1\subseteq H_0$, \\So $T_0H_0\subseteq t_0H_0$ and
$T_{1}^{-1}t_0H_0\subseteq H_0$, \\it means that

$t_0^{-1}T_0\subseteq H_0=H$  and $T_1^{-1}t_0\subseteq H_0=H$  for
some $t_0\in T_0$.\\Sufficiency, set $A=H\times \{0\}\cup
(t_0H)\times \{1\}$, \\because $<S_0\cup \{1_G\}\setminus
\{s_0\}>\leq H$, we can get $|S_0|-r_0=1$. \\Similarly,\\ because
$G_1=<S_1>\leq t_0Ht_0^{-1},T_1^{-1}t_0\subseteq H$ and
$t_0^{-1}T_0\subseteq H$, we can get that

$|S_1|-r_1=0$, $|T_0|-p=0$ and  $|T_1|-q=0$.\\ So $\lambda(X)\leq
|\omega^{+}(A)|=|H|<\delta(X)$. $\Box$
\\\\Analogously, we can achieve the following proposition 3.9, 3.10  and  3.11
easily.\\
\\\textbf{Proposition 3.9} $X=MD(G,S_0,S_1,T_0,T_1)$ is a strongly connected mixed
Cayley digraph, $X$ is Class 2 or  Class $4^{'}$ if and only if
\\(1) There exists a non-empty proper subgroup $H$ of $G$ such
that

$G_0=<S_0>\leq H$ and $|H|<\delta(X)$, and \\(2) There are two
distinct elements $t_0,t_0^{'}\in T_0$ such that

$G_1=<S_1>\leq t_0Ht_0^{-1}$, $T_1^{-1}t_0\subseteq H$,

$t_{0}^{'}H\cap t_0H=\varnothing$ and
$t_{0}^{-1}(T_{0}\setminus\{t_0^{'}\})\subseteq H$. $\Box$\\
\\\textbf{Proposition 3.10} $X=MD(G,S_0,S_1,T_0,T_1)$ is a strongly connected mixed
Cayley digraph, $X$ is Class 3 or  Class $3^{'}$ if and only if
\\(1) There exists a non-empty proper subgroup $H$ of $G$ and  some element $s_1\in
S_1$ such that

$<S_1\cup\{1_G\}\setminus\{s_1\}>\leq H$  and $|H|<\delta(X)$, and
\\(2) There is an element $t_1\in T_1$ such that

$G_0=<S_0>\leq t_1^{-1}Ht_1$, $T_0t_1^{-1}\subseteq H$ and
$t_1T_1^{-1}\subseteq H$. $\Box$\\
\\\textbf{Proposition 3.11} $X=MD(G,S_0,S_1,T_0,T_1)$ is a strongly connected mixed
Cayley digraph, $X$ is Class 4 or  Class $2^{'}$ if and only if
\\(1) There exists a non-empty proper subgroup $H$ of $G$ such
that

$G_1=<S_1>\leq H$ and $|H|<\delta(X)$, and\\(2) There are two
distinct elements $t_1\in T_1$ and $t_1^{'}\in T_1$ such that

$G_0=<S_0>\leq t_1^{-1}Ht_1$, $T_0t_1^{-1}\subseteq H$,

${t_{1}^{'}}^{-1}H \cap t_1^{-1}H=\varnothing$ and
$t_1(T_1^{-1}\backslash \{{t_{1}^{'}}^{-1}\}) \subseteq H$.
$\Box$\\\\Put the above
propositions together, we get the following theorem.\\
\\\textbf{Theorem 3.12} Let $X=MD(G,S_0,S_1,T_0,T_1)$ be a strongly
connected mixed Cayley digraph. Then $X$ is not max$-\lambda$ if and
only if $X$ satisfies one of the following
conditions:\\\textbf{\texttt{Condition 1.}} \\(1.1) There exists a
non-empty proper subgroup $H$ of $G$ and $S_0$ contains

an element $s_0$ such that

 $<S_0\cup \{1_G\}\setminus \{s_0\}>\leq
H$ and $|H|<\delta(X)$, and \\(1.2) There is an element $t_0\in T_0$
such that

$G_1=<S_1>\leq t_0Ht_0^{-1},$ $T_1^{-1}t_0\subseteq H$ and
$t_0^{-1}T_0\subseteq H$.\\\texttt{Condition 2.} \\(2.1) There
exists a non-empty proper subgroup $H$ of $G$ such that

$G_0=<S_0>\leq H$ and $|H|<\delta(X)$, and \\(2.2) There are two
distinct elements $t_0,$ $t_0^{'}\in T_0$ such that

$G_1=<S_1>\leq t_0Ht_0^{-1}$, $T_1^{-1}t_0\subseteq H$,

$t_{0}^{'}H\cap t_0H=\varnothing$ and
$t_{0}^{-1}(T_{0}\setminus\{t_0^{'}\})\subseteq
H$.\\\texttt{Condition 3.} \\(3.1) There exists a non-empty proper
subgroup $H$ of $G$ and  some element $s_1\in S_1$ such that

$<S_1\cup\{1_G\}\setminus\{s_1\}>\leq H$ and $|H|<\delta(X)$, and
\\(3.2) There is an element $t_1\in T_1$ such that

$G_0=<S_0>\leq t_1^{-1}Ht_1$, $T_0t_1^{-1}\subseteq H$ and
$t_1T_1^{-1}\subseteq H$.\\\texttt{Condition 4.} \\(4.1) There
exists a non-empty proper subgroup $H$ of $G$ such that

$G_1=<S_1>\leq H$ and $|H|<\delta(X)$, and \\(4.2) There are two
distinct elements $t_1\in T_1$ and $t_1^{'}\in T_1$ such that

$G_0=<S_0>\leq t_1^{-1}Ht_1$, $T_0t_1^{-1}\subseteq H,$

${t_{1}^{'}}^{-1}H \cap t_1^{-1}H=\varnothing$ and
$t_1(T_1^{-1}\backslash
\{{t_{1}^{'}}^{-1}\}) \subseteq H$. $\Box$\\

So all the strongly connected mixed Cayley digraphs but a few
exceptions are max-$\lambda$ .
\section{Super arc-connectivity of mixed Cayley digraph} If a digraph isn$^{'}$t max-$\lambda$, it is also not super-$\lambda$. So in this section, we
prove that the mixed Cayley digraph, which is max-$\lambda$ but not
super-$\lambda$, is only a few exceptions. A \emph{weak path} of a
digraph $X$ is a sequence $u_0,...,u_r$ of distinct vertices such
that for $i=1,...,r,$ either $(u_{i-1},u_i)$ or $(u_i,u_{i-1})$ is
an arc of $X$. A directed graph is \emph{weakly connected} if any
two vertices can be joined by a weak path. The following proposition
is clearly true.\\\\\textbf{proposition 4.1} If
$X=MD(G,S_0,S_1,T_0,T_1)$ is max$-\lambda$, but not super$-\lambda$,
let A be a $\lambda-$super atom, then $Y=X[A]$ is weakly connected.
$\Box$\\
\\\textbf{Lemma 4.2} If $X=MD(G,S_0,S_1,T_0,T_1)$ is
max$-\lambda$, but not super$-\lambda$, let A be a
$\lambda-$superatom, If $\delta(X)=1$, $A\cap
X_i\neq\varnothing$.\\\\Proof. By contradiction. Because $X$ is
strongly connected, $|T_0|\geq 1$ and $|T_1|\geq 1$. If
$\delta(X)=1$ , one of the following conditions holds:

 $(1)$
$|S_0|=0$, $|T_0|=1$ or $|T_1|=1$,

 (2) $|S_1|=0$, $|T_0|=1$ or
$|T_1|=1$.\\ Without loss generality, suppose (1) holds, then
$A\nsubseteq X_0$, thus $\lambda(X)\geq |A||T_1|\geq|A|\geq2$, a
contradiction. $\Box$\\\\\textbf{Lemma 4.3} Let
$X=MD(G,S_0,S_1,T_0,T_1)$ is max$-\lambda$, but not super$-\lambda$,
let $A$ be a $\lambda-$superatom, then $|A|\geq
\delta(X)$.\\\\Proof. Suppose $A$ is a positive $\lambda-$superatom.
Then

$\lambda(X)=|\omega_X^{+}(A)|\geq|A|(\delta-(|A|-1))=|A|(\delta-|A|+1)$,
\\we can verify that $\lambda(X)>\delta(X)$ when
$2\leq|A|<\delta(X)$,
a contradiction.So $|A|\geq\delta(X)$. $\Box$\\

Let $X=MD(G,S_0,S_1,T_0,T_1)=Cayley(G\times\{0\},S_0)\cup
Cayley(G\times\{1\},S_1)\cup BD(G,T_0,T_1) $ be a strongly connected
mixed digraph. There is a class of special mixed Cayley digraph,
which is that one of $Cayley(G\times\{0\},S_0)$ and $
Cayley(G\times\{1\},S_1)$ is a union of disjoint directed cycles,
and the other is a union of disjoint directed cycles with length
two, and $BD(G,T_0,T_1)$ is a union of disjoint directed cycles with
length two. The class of special mixed Cayley
digraphs is denoted by $\mathcal {F}$.\\
\\\textbf{Lemma 4.4} Let $X=MD(G,S_0,S_1,T_0,T_1)$ is
max$-\lambda$ but not super$-\lambda$, if $X$ is neither a directed
cycle nor a cycle and $X$ doesn$^{'}$t belong to $\mathcal {F}$,
then distinct positive(respectively, negative) $\lambda-$superatoms
of $X$ are vertex disjoint.\\\\Proof. We suppose $\delta(X)\leq2$,
since  the lemma is true when $\delta(X)\geq3$ by theorem 2.4.

Let $A$ and $B$ be two distinct positive $\lambda-$superatoms.  If
$A\cap B\neq\varnothing$, by proposition 2.1, $A\cap B,A\cup
B,A\setminus B,B\setminus A$ are arc fragments of $X$. Because each
of $A\cap B,A\setminus B$ and $B\setminus A$ is a proper subset of a
$\lambda-$superatom, we achieve that

$|A\cap B|=1$, $|A\setminus B|=1$ and $|B\setminus A|=1$. \\So
assume $A=\{u,v\}$, $B=\{v,w\}$ with $u\neq w$, thus

 $d_{X[A]}^+(u)=d_{X[A]}^-(v)\leq1$,
$d_{X[A]}^-(u)=d_{X[A]}^+(v)\leq1$,

$d_{X[B]}^+(v)=d_{X[B]}^-(w)\leq1$ and
$d_{X[B]}^-(v)=d_{X[B]}^+(w)\leq1$.\\\textbf{Case 1} $A$,
$B\subseteq X_0$ or $A$, $B\subseteq X_1$.

$2\geq\lambda(X)=|\omega_X^+(A\cup B)|\geq3min(|T_0|,|T_1|)\geq3$, a
contradiction.\\\textbf{Case 2} $A\cap X_i\neq\varnothing$ and
$B\cap X_i\neq\varnothing$.

When $\delta(X)=1$. \\Since $A\cap B,A\setminus B$ and $B\setminus
A$ are arc-fragments of $X$, $d_X^+(u)=d_X^+(v)=d_X^+(w)=1$, so
$|T_0|=1$, $|T_1|=1$ and $|S_0|=|S_1|=0$. And because $X$ is a
strongly connected digraph, we can get $X$ is a directed cycle, a
contradiction.

When $\delta(X)=2$, then

$d_X^+(u)=d_X^+(v)=d_X^+(w)=2$. \\Because $X[A]$ and $X[B]$ are
weakly connected and $A$, $B$ and $A\cup B$ are arc fragments, we
can deduce

$|T_0|=|T_1|=2$, $T_0=T_1$ and $|S_0|=|S_1|=0$.\\ Because $X$ is
strongly connected, $X$ is a cycle, a contradiction.\\\textbf{Case
3} $A\cap X_i\neq\varnothing$ and either $B\subseteq X_0$ or
$B\subseteq X_1$ ( $B\cap X_i\neq\varnothing$ and $A\subseteq X_0$
or $A\subseteq X_1$). \\By lemma 4.2, we can get

$\delta(X)\geq2$, so $\delta(X)=2$. \\Since $A$, $B$ and $A\cup B$
are arc fragment, We can get that

$|T_0|=|T_1|=|S_0|=|S_1|=1$ such that $S_1^{-1}=S_1$ or
$S_0^{-1}=S_0 $ and $T_0=T_1$, \\thus $X$ belongs to $\mathcal {F}$,
a contradiction. $\Box$\\\\\textbf{Lemma 4.5} Let
$X=MD(G,S_0,S_1,T_0,T_1)$ is max$-\lambda$ but not super$-\lambda$,
if $X$ is neither a directed cycle nor a cycle and $X$ doesn$^{'}$t
belong to $\mathcal {F}$. Let $A$ be a $\lambda-$superatom of $X$,
then\\ (1) If $A\subseteq X_i$, let $A=H\times
\{i\},i=0,1,H\subseteq G$. And let $Y=X[A]$ be the subgraph of $X$
induced by $A$, then

(i) \emph{Aut(Y)} acts transitively on $A$, and

(ii) If $A$ contains $(1_G, 0)$ or $(1_G, 1)$, then $H$ is a
subgroup of $G$.\\(2) If $A_i=A\cap
X_i=H_i\times\{i\}\neq\varnothing$ where $H_i\subseteq G$, let
$Y_i=X[A_i]$ be the subgraphs of $X$ induced by$A_i$, then

(i) \emph{Aut($Y_i$)} acts transitively on $A_i$ for $i=0,1$, and

(ii) If $A_i$ contains $(1_G, i) (i=0, 1)$,then $H_i$ is a subgroup
of $G$.\\\\Proof. (1) Without loss of generality,suppose $A\subseteq
X_0$, then $X[A]$ is the subdigraph of $Cayley(G\times\{0\},S_0)$
induced by $H\times\{0\}\subseteq G\times\{0\}$, where
$A=H\times\{0\}$. \\By lemma 4.4, $A$ is an imprimitive block of
$Cayley(G\times\{0\},S_0)$. \\So by theorem 2.3, we have that (1)
holds.\\(2) The proof is similar to the proof of lemma 3.3.
$\Box$\\\\\textbf{Lemma 4.6} Let $X=MD(G,S_0,S_1,T_0,T_1)$ be
max$-\lambda$ but not super$-\lambda$, if $X$ is neither a directed
cycle nor a cycle and $X$ doesn$^{'}$t belong to $\mathcal {F}$. Let
$A$ be a $\lambda-$super atom of $X$ with $A_i=A\cap X_i\neq
\varnothing$,
 then\\(1) $V(X)$ is a disjoint union of distinct positive(negative)
$\lambda-$super atoms, and\\(2) $|A_0|=|A_1|$. $\Box$\\

The proof of lemma 4.6 is similar to the proof of lemma
3.4.\\\\\textbf{proposition 4.7} Let $X=MD(G,S_0,S_1,T_0,T_1)$ be
max$-\lambda$ but not super$-\lambda$, if $X$ is neither a directed
cycle nor a cycle and $X$ doesn$^{'}$t belong to $\mathcal {F}$. Let
$A$ be a $\lambda-$super atom of $X$ with $A_i=A\cap X_i\neq
\varnothing$. Set $A_{i}=A\cap X_{i}=H_{i}\times\{i\}$, for $i=0,1$,
where $H_{i}\subseteq G$, then\\$(1)$ If $(1_{G},0)\in A_{0}$, then
$H_{1}=t_{0}H_{0}$ for some $t_{0}\in T_{0}$,
furthermore,

$X_{0}=\cup_{i=1}^{k}(H_{0}g_{i})\times\{0\}$,

$X_{1}=\cup_{i=1}^{k}(t_{0}H_{0}g_{i})\times\{1\}$, \\where
$R(g_{i})(A)\cap R(g_{j})(A)\neq\varnothing$ if and only if $i=j$
for $1\leq i,j \leq k$.\\(2) If $(1_{G},1)\in A_{1}$, then
$H_{0}=t_{1}^{-1}H_{1}$ for some $t_{1}\in T_{1}$, furthermore,

$X_{0}=\cup_{i=1}^{k}(t_{1}^{-1}H_{1}g_{i})\times\{0\}$,

$X_{1}=\cup_{i=1}^{k}(H_{1}g_{i})\times\{1\}$,\\ where
$R(g_{i})(A)\cap R(g_{j})(A)\neq\varnothing$ if and only if
$i=j$ for $1\leq i,j \leq k$. $\Box$\\\\
The proof is similar to the lemma 3.5.\\

We give two classes of digraphs which aren$^{'}$t super-$\lambda$.
The first of class digraphs consists of the strongly connected mixed
Cayley digraphs $X=MD(G,S_0,S_1,T_0,T_1)$ which contain
$\lambda-$superatoms lining in $X_0$ or $X_1$. This class of
digraphs is denoted by $\mathcal {G}$. The second class of digraphs
consists of the strongly connected mixed Cayley digraph
$X=MD(G,S_0,S_1,T_0,T_1)$ all of whose $\lambda-$superatoms contain
at least one vertex of $X_0$ and $X_1$ respectively, denoted by
$\mathcal {L}$.
\\\\\textbf{Theorem 4.8} Let
$X=MD(G,S_0,S_1,T_0,T_1)$ be max$-\lambda$, but not super$-\lambda$,
if $X$ is neither a directed cycle nor a cycle and $X$ doesn$^{'}$t
belong to $\mathcal {F}$, then $X$ belongs to $\mathcal {G}$ if and
only if $X$ satisfies one of the following conditions:\\(1)
$|T_0|=1$ or $|T_1|=1$,$1\leq|S_0|\leq|S_1|$ and $S_0\cup\{1_G\}\leq
G$.\\(2) $|T_0|=1$ or $|T_1|=1$,$1\leq|S_1|\leq|S_0|$ and
$S_1\cup\{1_G\}\leq G$.\\\\Proof. Necessity. \\Because $A\subseteq
X_0$ or $A\subseteq X_1$, we have $\delta(X)\geq2$ by the lemma
4.2.\\\textbf{1.1} $A$ is a positive $\lambda-$superatom.

If $A\subseteq X_0$ and $(1_G,0)\in A$, then $A=H\times \{0\}$ and
$H\leq G$. \\By lemma 4.5, $Y=X[A]$ is a regular digraph, so we set
$Y$ is $r_0$ regular digraph, then

$\lambda(X)=|\omega_X^{+}(A)|=|A|(|S_0|+|T_0|-r_0)\geq
\delta(X)(|S_0|+|T_0|-r_0)$. \\Since $|T_0|\geq1$ and
$\lambda(X)=\delta(X)$, we have that

$|A|=\delta(X),|S_0|-r_0=0$ and $|T_0|=1$. \\So $Y$ is a $|S_0|$
regular digraph with order $\delta(X)$. Thus

$|S_0|\leq\delta(X)-1=min\{|S_0|+1, |S_1|+1, |S_0|+|T_1|,
|S_1|+|T_1|\}-1\leq min\{|S_0|+1, |S_1|+1\}-1$, \\so we can get
$|S_0|\leq|S_1|$ and $Y\cong K_{\delta(X)}$. So if $X$ is not
super$-\lambda$ then

$|T_0|=1,1\leq|S_0|\leq|S_1|$ and $S_0\cup\{1_G\}\leq G$.
\\Similarly, if $A\subseteq X_1$, we can prove that if $X$ is not super$-\lambda$, then

$|T_1|=1,1\leq|S_1|\leq|S_0|$ and $S_1\cup\{1_G\}\leq
G$.\\\textbf{1.2}  $A$ is a negative $\lambda-$super atom

If $A\subseteq X_0$ and $(1_G,0)\in A$. Let $A=H\times\{0\}$, then
$H\leq G$, by lemma 4.5, $Y=X[A]$ is a regular digraph. We set $Y$
is $r_0$ regular digraph, then

$\lambda(X)=|\omega_X^{-1}(A)|=|A|(|S_0|+|T_1|-r_0)\geq\delta(X)(|S_0|+|T_1|-r_0)$.\\Since
$|T_1|\geq1$ and $\lambda(X)=\delta(X)$, we have that

$|A|=\delta(X), |T_1|=1$ and $|S_0|=r_0$.\\ So $Y$ is a
$|S_0|-$regular digraph with order $\delta(X)$, thus

$|S_0|\leq\delta(X)-1\leq min\{|S_0|+1,|S_1|+1\}-1$. \\So we have
that

$|S_0|\leq|S_1|$ and  $|S_0|=\delta(X)-1\geq1$. \\So $Y\cong
K_{\delta(X)}$  and  $H=S_0\cup\{1_G\}\leq G$. \\So if $X$ is not
super$-\lambda$ then

$|T_1|=1,1\leq|S_0|\leq|S_1|$ and $S_0\cup\{1_G\}\leq G$.\\
Similarly, if $A\subseteq X_1$, we can achieve that if $X$ is not
super$-\lambda$ then

$|T_0|=1, 1\leq|S_1|\leq|S_0|$ and $S_1\cup \{1_g\}\leq
G$.\\Sufficiency. For condition (1), because of $|T_0|\geq 1$ and
$|T_1|\geq 1$, we have

$\delta(X)=|S_0|+1$. \\Set $A=S_0\times \{0\}\cup \{(1_G,0)\}$, then

min\{$\omega^{+}(A), \omega^{-}(A)$\}=$|A|\times$min$\{|T_0|,
|T_1|\}$=$|A|$=$|S_0|+1$=$\delta(X)$, \\and because

$|A|=|S_0\times \{0\}\cup \{(1_G,0)\}|\geq 2$, \\so $A$ is a
nontrivial $\lambda-$fragment.\\ Condition (2) is similar to
condition
(1). $\Box$\\

For the class of $\mathcal {L}$, by lemma 4.5 and proposition 1.5,
we can prove that \emph{Aut(Y)} acts transitively both on $A_i$
where $Y=X[A]$. Thus if we set $Y^{'}=Y\setminus\{E(Y_0)\cup
E(Y_1)\}$ where $Y_i=X[A_i]$, we can easily prove that

$d_{Y^{'}}^{+}((g,0))=d_{Y^{'}}^{-}((g,1))$ and
$d_{Y^{'}}^{-}((g,0))=d_{Y^{'}}^{+}((g,1))$.\\ So we set

$d_{Y^{'}}^{+}((g,0))=d_{Y^{'}}^{-}((g,1))=p$ and

$d_{Y^{'}}^{-}((g,0))=d_{Y^{'}}^{+}((g,1))=q$, and

let $Y_i$ is $r_i$-regular digraph  for $i=0,1$ .\\There are some special Classes of digraphs of $\mathcal{L}$ as follows:

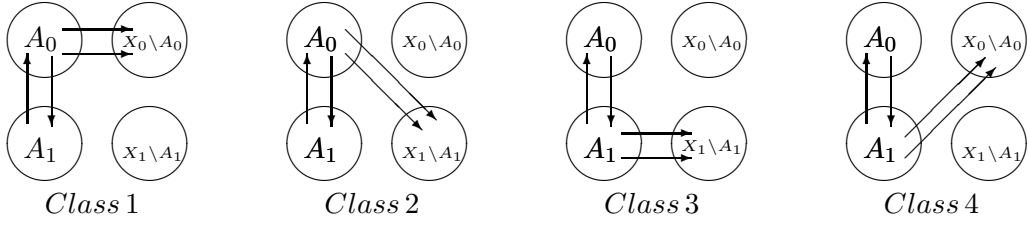
\begin{figure}[h,t]
\setlength{\unitlength}{0.46mm}
\begin{center}
\begin{picture}(290,56)

\put(0,10){\circle{20}}
\put(0,40){\circle{20}}

\put(30,10){\circle{20}}
\put(30,40){\circle{20}}

\put(80,10){\circle{20}}
\put(80,40){\circle{20}}

\put(110,10){\circle{20}}
\put(110,40){\circle{20}}

\put(160,10){\circle{20}}
\put(160,40){\circle{20}}

\put(190,10){\circle{20}}
\put(190,40){\circle{20}}

\put(240,10){\circle{20}}
\put(240,40){\circle{20}}

\put(270,10){\circle{20}}
\put(270,40){\circle{20}}

\put(0,-10){$Class\,1$}
\put(80,-10){$Class\,2$}
\put(160,-10){$Class\,3$}
\put(240,-10){$Class\,4$}
\put(-6,6){$A_1$}
\put(-6,38){$A_0$}

\put(74,6){$A_1$}
\put(74,38){$A_0$}

\put(154,6){$A_1$}
\put(154,38){$A_0$}

\put(234,6){$A_1$}
\put(234,38){$A_0$}

\put(74,6){$A_1$}
\put(74,38){$A_0$}

\put(154,6){$A_1$}
\put(154,38){$A_0$}

\put(234,6){$A_1$}
\put(234,38){$A_0$}
\put(22,6){\tiny$X_1\backslash A_1$}
\put(22,38){\tiny $X_0\backslash A_0$}

\put(102,6){\tiny$X_1\backslash A_1$}
\put(102,38){\tiny $X_0\backslash A_0$}

\put(182,8){\tiny$X_1\backslash A_1$}
\put(182,38){\tiny $X_0\backslash A_0$}

\put(262,6){\tiny$X_1\backslash A_1$}
\put(262,38){\tiny $X_0\backslash A_0$}
\put(-5,16){\vector(0,1){20}}
\put(2,35){\vector(0,-1){20}}

\put(75,16){\vector(0,1){20}}
\put(82,35){\vector(0,-1){20}}

\put(155,16){\vector(0,1){20}}
\put(162,35){\vector(0,-1){20}}

\put(235,16){\vector(0,1){20}}
\put(242,35){\vector(0,-1){20}}
\put(5,36){\vector(1,0){20}}
\put(5,43){\vector(1,0){20}}

\put(165,6){\vector(1,0){20}}
\put(165,13){\vector(1,0){20}}

\put(86,36){\vector(1,-1){22}}
\put(86,43){\vector(1,-1){26}}
\put(246,6){\vector(1,1){26}}
\put(246,12){\vector(1,1){23}}
\end{picture}
\end{center}
\end{figure}

Where $A$ is a $\lambda-$superatom and $A_i=X_i\cap A$, and\\ Class
1 satisfies that $|S_0|-r_0=1$, $|S_1|-r_1=0$, $|T_0|-p=0$ and
$|T_1|-q=0$, and\\ Class 2 satisfies that $|S_0|-r_0=0$,
$|S_1|-r_1=0$, $|T_0|-p=1$ and $|T_1|-q=0$, and \\Class 3 satisfies
that $|S_0|-r_0=0$, $|S_1|-r_1=1$, $|T_0|-p=0$ and $|T_1|-q=0$,
and\\
Class 4 satisfies that $|S_0|-r_0=0$, $|S_1|-r_1=0$, $|T_0|-p=0$
and $|T_1|-q=1$, and\\ all of the above digraphs satisfy that
$|A_0|=|A_1|=\delta(X).$

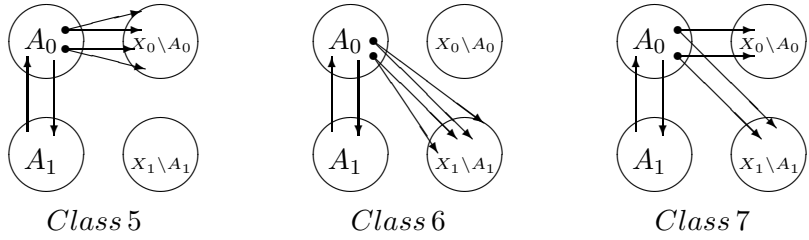
\begin{figure}[h,t]
\setlength{\unitlength}{0.5mm}
\begin{center}
\begin{picture}(180,56)

\put(0,10){\circle{20}}
\put(0,40){\circle{20}}

\put(30,10){\circle{20}}
\put(30,40){\circle{20}}

\put(80,10){\circle{20}}
\put(80,40){\circle{20}}

\put(110,10){\circle{20}}
\put(110,40){\circle{20}}

\put(160,10){\circle{20}}
\put(160,40){\circle{20}}

\put(190,10){\circle{20}}
\put(190,40){\circle{20}}

\put(0,-10){$Class\,5$}
\put(80,-10){$Class\,6$}
\put(160,-10){$Class\,7$}
\put(-6,6){$A_1$}
\put(-6,38){$A_0$}

\put(74,6){$A_1$}
\put(74,38){$A_0$}

\put(154,6){$A_1$}
\put(154,38){$A_0$}
\put(22,6){\tiny$X_1\backslash A_1$}
\put(22,38){\tiny $X_0\backslash A_0$}

\put(102,6){\tiny$X_1\backslash A_1$}
\put(102,38){\tiny $X_0\backslash A_0$}

\put(182,6){\tiny$X_1\backslash A_1$}
\put(182,38){\tiny $X_0\backslash A_0$}
\put(-5,16){\vector(0,1){20}}
\put(2,35){\vector(0,-1){20}}

\put(75,16){\vector(0,1){20}}
\put(82,35){\vector(0,-1){20}}

\put(155,16){\vector(0,1){20}}
\put(162,35){\vector(0,-1){20}}

\put(5,38){\vector(4,-1){21}}
\put(5,38){\vector(1,0){18}}
\put(5,38){\circle*{2}}
\put(5,43){\vector(1,0){20}}
\put(5,43){\vector(4,1){20}}
\put(5,43){\circle*{2}}

\put(86,36){\vector(1,-1){22}}
\put(86,36){\vector(2,-3){17}}
\put(86,36){\circle*{2}}
\put(86,40){\vector(1,-1){26}}
\put(86,40){\vector(4,-3){29}}
\put(86,40){\circle*{2}}

\put(166,36){\vector(1,-1){22}}
\put(166,43){\vector(1,-1){26}}
\put(166,36){\vector(1,0){20}}
\put(166,43){\vector(1,0){20}}
\put(166,36){\circle*{2}}
\put(166,43){\circle*{2}}
\end{picture}
\end{center}
\end{figure}

Where $A$ is a $\lambda-$superatom and $A_i=X_i\cap A$, and\\ Class
5 satisfies that $|S_0|-r_0=2$, $|S_1|-r_1=0$, $|T_0|-p=0$ and
$|T_1|-q=0$, and\\ Class 6 satisfies that $|S_0|-r_0=0$,
$|S_1|-r_1=0$, $|T_0|-p=2$ and $|T_1|-q=0$, and \\Class 7 satisfies
that $|S_0|-r_0=1$, $|S_1|-r_1=0$, $|T_0|-p=1$ and $|T_1|-q=0$,
and\\
 all of the above digraphs satisfy that $|A_0|=|A_1|=\delta(X)/2.$
\begin{figure}[h,t]
\setlength{\unitlength}{0.5mm}
\begin{center}
\begin{picture}(180,56)

\put(0,10){\circle{20}}
\put(0,40){\circle{20}}

\put(30,10){\circle{20}}
\put(30,40){\circle{20}}

\put(80,10){\circle{20}}
\put(80,40){\circle{20}}

\put(110,10){\circle{20}}
\put(110,40){\circle{20}}

\put(160,10){\circle{20}}
\put(160,40){\circle{20}}

\put(190,10){\circle{20}}
\put(190,40){\circle{20}}

\put(0,-10){$Class\,8$}
\put(80,-10){$Class\,9$}
\put(160,-10){$Class\,10$}
\put(-6,6){$A_1$}
\put(-6,38){$A_0$}

\put(74,6){$A_1$}
\put(74,38){$A_0$}

\put(154,6){$A_1$}
\put(154,38){$A_0$}
\put(22,9){\tiny$X_1\backslash A_1$}
\put(22,38){\tiny $X_0\backslash A_0$}

\put(102,9){\tiny$X_1\backslash A_1$}
\put(102,38){\tiny $X_0\backslash A_0$}

\put(182,9){\tiny$X_1\backslash A_1$}
\put(182,38){\tiny $X_0\backslash A_0$}
\put(-5,16){\vector(0,1){20}}
\put(2,35){\vector(0,-1){20}}

\put(75,16){\vector(0,1){20}}
\put(82,35){\vector(0,-1){20}}

\put(155,16){\vector(0,1){20}}
\put(162,35){\vector(0,-1){20}}

\put(5,8){\vector(4,-1){21}}
\put(5,8){\vector(1,0){18}}
\put(5,8){\circle*{2}}
\put(5,13){\vector(1,0){20}}
\put(5,12){\vector(4,1){20}}
\put(5,12){\circle*{2}}

\put(86,11){\vector(1,1){22}}
\put(86,11){\vector(4,3){28}}
\put(86,11){\circle*{2}}
\put(86,15){\vector(2,3){16}}
\put(86,15){\vector(1,1){20}}
\put(86,15){\circle*{2}}

\put(166,6){\vector(1,1){28}}
\put(166,6){\vector(1,0){18}}
\put(166,6){\circle*{2}}

\put(166,13){\vector(1,1){23}}
\put(166,13){\vector(1,0){18}}
\put(166,13){\circle*{2}}
\end{picture}
\end{center}
\end{figure}

Where $A$ is a $\lambda-$superatom and $A_i=X_i\cap A$, and\\ Class
8 satisfies that $|S_0|-r_0=0$, $|S_1|-r_1=2$, $|T_0|-p=0$ and
$|T_1|-q=0$, and\\ Class 9 satisfies that $|S_0|-r_0=0$,
$|S_1|-r_1=0$, $|T_0|-p=0$ and $|T_1|-q=2$, and \\Class 10
satisfies that $|S_0|-r_0=1$, $|S_1|-r_1=0$, $|T_0|-p=0$ and
$|T_1|-q=1$,
and\\
 all of the above digraphs satisfy that $|A_0|=|A_1|=\delta(X)/2.$\\\\\\

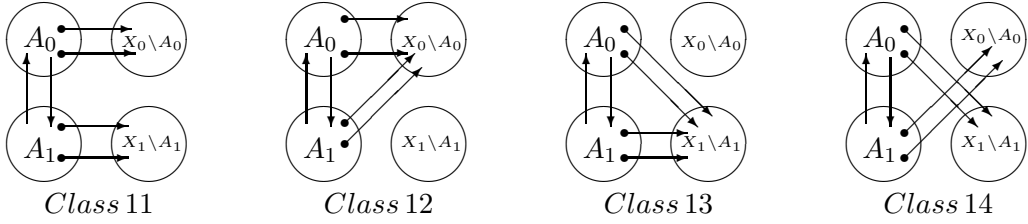
\begin{figure}[h,t]
\setlength{\unitlength}{0.46mm}
\begin{center}
\begin{picture}(260,56)

\put(0,10){\circle{20}}
\put(0,40){\circle{20}}

\put(30,10){\circle{20}}
\put(30,40){\circle{20}}

\put(80,10){\circle{20}}
\put(80,40){\circle{20}}

\put(110,10){\circle{20}}
\put(110,40){\circle{20}}

\put(160,10){\circle{20}}
\put(160,40){\circle{20}}

\put(190,10){\circle{20}}
\put(190,40){\circle{20}}

\put(240,10){\circle{20}}
\put(240,40){\circle{20}}

\put(270,10){\circle{20}}
\put(270,40){\circle{20}}

\put(0,-10){$Class\,11$}
\put(80,-10){$Class\,12$}
\put(160,-10){$Class\,13$}
\put(240,-10){$Class\,14$}
\put(-6,6){$A_1$}
\put(-6,38){$A_0$}

\put(74,6){$A_1$}
\put(74,38){$A_0$}

\put(154,6){$A_1$}
\put(154,38){$A_0$}
\put(234,6){$A_1$}
\put(234,38){$A_0$}
\put(22,9){\tiny$X_1\backslash A_1$}
\put(22,38){\tiny $X_0\backslash A_0$}

\put(102,9){\tiny$X_1\backslash A_1$}
\put(102,38){\tiny $X_0\backslash A_0$}

\put(182,9){\tiny$X_1\backslash A_1$}
\put(182,38){\tiny $X_0\backslash A_0$}

\put(262,9){\tiny$X_1\backslash A_1$}
\put(262,40){\tiny $X_0\backslash A_0$}
\put(-5,16){\vector(0,1){20}}
\put(2,35){\vector(0,-1){20}}

\put(75,16){\vector(0,1){20}}
\put(82,35){\vector(0,-1){20}}

\put(155,16){\vector(0,1){20}}
\put(162,35){\vector(0,-1){20}}

\put(235,16){\vector(0,1){20}}
\put(242,35){\vector(0,-1){20}}
\put(5,36){\vector(1,0){21}}
\put(5,43){\vector(1,0){20}}
\put(5,36){\circle*{2}}
\put(5,43){\circle*{2}}

\put(5,6){\vector(1,0){19}}
\put(5,6){\circle*{2}}
\put(5,15){\vector(1,0){19}}
\put(5,15){\circle*{2}}

\put(86,36){\vector(1,0){18}}
\put(86,46){\vector(1,0){18}}
\put(86,36){\circle*{2}}
\put(86,46){\circle*{2}}
\put(86,10){\vector(1,1){22}}
\put(86,16){\vector(1,1){20}}
\put(86,10){\circle*{2}}
\put(86,16){\circle*{2}}

\put(166,36){\vector(1,-1){21}}
\put(166,6){\vector(1,0){18}}
\put(166,6){\circle*{2}}
\put(166,36){\circle*{2}}

\put(166,43){\vector(1,-1){25}}
\put(166,13){\vector(1,0){18}}
\put(166,13){\circle*{2}}
\put(166,43){\circle*{2}}

\put(246,6){\vector(1,1){28}}
\put(246,6){\circle*{2}}
\put(246,13){\vector(1,1){25}}
\put(246,13){\circle*{2}}

\put(246,36){\vector(1,-1){21}}
\put(246,43){\vector(1,-1){25}}
\put(246,36){\circle*{2}}
\put(246,43){\circle*{2}}
\end{picture}
\end{center}
\end{figure}

Where $A$ is a $\lambda-$superatom and $A_i=X_i\cap A$, and\\ Class
11 satisfies that $|S_0|-r_0=1$, $|S_1|-r_1=1$, $|T_0|-p=0$ and
$|T_1|-q=0$, and\\ Class 12 satisfies that $|S_0|-r_0=1$,
$|S_1|-r_1=0$, $|T_0|-p=1$ and $|T_1|-q=1$, and \\Class 13
satisfies that $|S_0|-r_0=0$, $|S_1|-r_1=1$, $|T_0|-p=1$ and
$|T_1|-q=0$,
and\\
Class 14 satisfies that $|S_0|-r_0=0$, $|S_1|-r_1=0$, $|T_0|-p=1$
and $|T_1|-q=1$, and\\ all of the above digraphs satisfy that
$|A_0|=|A_1|=\delta(X)/2.$

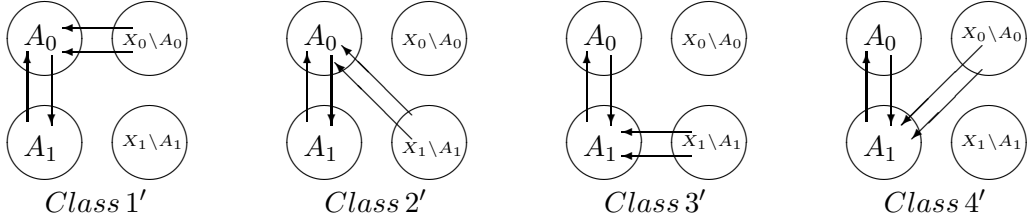
\begin{figure}[h,t]
\setlength{\unitlength}{0.46mm}
\begin{center}
\begin{picture}(290,56)

\put(0,10){\circle{20}}
\put(0,40){\circle{20}}

\put(30,10){\circle{20}}
\put(30,40){\circle{20}}

\put(80,10){\circle{20}}
\put(80,40){\circle{20}}

\put(110,10){\circle{20}}
\put(110,40){\circle{20}}

\put(160,10){\circle{20}}
\put(160,40){\circle{20}}

\put(190,10){\circle{20}}
\put(190,40){\circle{20}}

\put(240,10){\circle{20}}
\put(240,40){\circle{20}}

\put(270,10){\circle{20}}
\put(270,40){\circle{20}}

\put(0,-10){$Class\,1'$}
\put(80,-10){$Class\,2'$}
\put(160,-10){$Class\,3'$}
\put(240,-10){$Class\,4'$}
\put(-6,6){$A_1$}
\put(-6,38){$A_0$}

\put(74,6){$A_1$}
\put(74,38){$A_0$}

\put(154,6){$A_1$}
\put(154,38){$A_0$}

\put(234,6){$A_1$}
\put(234,38){$A_0$}

\put(22,9){\tiny$X_1\backslash A_1$}
\put(22,38){\tiny $X_0\backslash A_0$}

\put(102,7){\tiny$X_1\backslash A_1$}
\put(102,38){\tiny $X_0\backslash A_0$}

\put(182,9){\tiny$X_1\backslash A_1$}
\put(182,38){\tiny $X_0\backslash A_0$}

\put(262,9){\tiny$X_1\backslash A_1$}
\put(262,40){\tiny $X_0\backslash A_0$}

\put(-5,16){\vector(0,1){20}}
\put(2,35){\vector(0,-1){20}}

\put(75,16){\vector(0,1){20}}
\put(82,35){\vector(0,-1){20}}

\put(155,16){\vector(0,1){20}}
\put(162,35){\vector(0,-1){20}}

\put(235,16){\vector(0,1){20}}
\put(242,35){\vector(0,-1){20}}
\put(25,36){\vector(-1,0){20}}
\put(25,43){\vector(-1,0){20}}

\put(185,6){\vector(-1,0){20}}
\put(185,13){\vector(-1,0){20}}

\put(105,11){\vector(-1,1){22}}
\put(105,18){\vector(-1,1){20}}
\put(268,31){\vector(-1,-1){20}}
\put(268,38){\vector(-1,-1){23}}
\end{picture}
\end{center}
\end{figure}

Where $A$ is a $\lambda-$superatom and $A_i=X_i\cap A$, and\\ Class
$1^{'}$ satisfies that $|S_0|-r_0=1$, $|S_1|-r_1=0$, $|T_0|-p=0$ and
$|T_1|-q=0$, and\\ Class $2^{'}$ satisfies that $|S_0|-r_0=0$,
$|S_1|-r_1=0$, $|T_0|-p=0$ and $|T_1|-q=1$, and \\Class $3^{'}$
satisfies that $|S_0|-r_0=0$, $|S_1|-r_1=1$, $|T_0|-p=0$ and
$|T_1|-q=0$,
and\\
Class $4^{'}$ satisfies that $|S_0|-r_0=0$, $|S_1|-r_1=0$,
$|T_0|-p=1$ and $|T_1|-q=0$, and\\ all of the above digraphs satisfy
that $|A_0|=|A_1|=\delta(X).$
Clearly,\\ Class $1^{'}$ is equivalent to Class 1,\\ Class $2^{'}$
is equivalent to Class 4, \\Class $3^{'}$ is equivalent to Class
3, and \\Class $4^{'}$ is equivalent to Class 2.

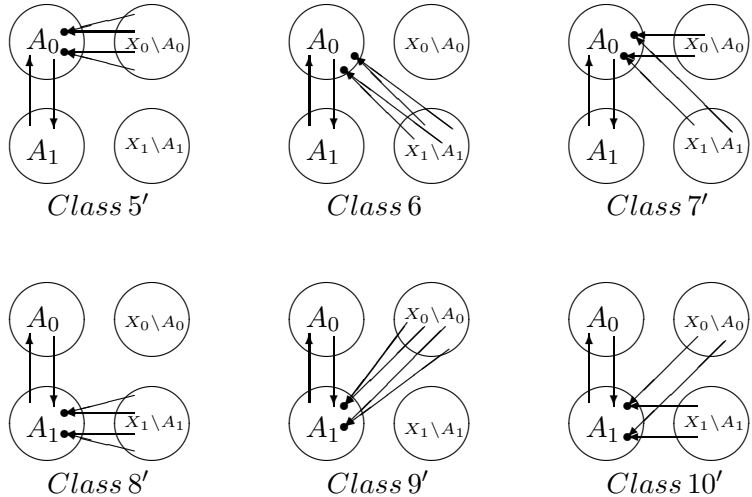
\begin{figure}[h,t]
\setlength{\unitlength}{0.46mm}
\begin{center}
\begin{picture}(200,136)

\put(0,10){\circle{20}}
\put(0,40){\circle{20}}

\put(30,10){\circle{20}}
\put(30,40){\circle{20}}

\put(80,10){\circle{20}}
\put(80,40){\circle{20}}

\put(110,10){\circle{20}}
\put(110,40){\circle{20}}

\put(160,10){\circle{20}}
\put(160,40){\circle{20}}

\put(190,10){\circle{20}}
\put(190,40){\circle{20}}
\put(-6,6){$A_1$}
\put(-6,38){$A_0$}

\put(74,6){$A_1$}
\put(74,38){$A_0$}

\put(154,6){$A_1$}
\put(154,38){$A_0$}

\put(22,9){\tiny$X_1\backslash A_1$}
\put(22,38){\tiny $X_0\backslash A_0$}

\put(102,7){\tiny$X_1\backslash A_1$}
\put(102,40){\tiny $X_0\backslash A_0$}

\put(182,9){\tiny$X_1\backslash A_1$}
\put(182,38){\tiny $X_0\backslash A_0$}

\put(-5,16){\vector(0,1){20}}
\put(2,35){\vector(0,-1){20}}

\put(75,16){\vector(0,1){20}}
\put(82,35){\vector(0,-1){20}}

\put(155,16){\vector(0,1){20}}
\put(162,35){\vector(0,-1){20}}

\put(5,7){\circle*{2}}
\put(5,13){\circle*{2}}
\put(25,7){\vector(-1,0){20}}
\put(25,2){\vector(-4,1){20}}
\put(25,13){\vector(-1,0){20}}
\put(25,18){\vector(-4,-1){20}}

\put(85,9){\circle*{2}}
\put(105,29){\vector(-1,-1){20}}
\put(105,29){\line(1,1){9}}
\put(105,24){\vector(-4,-3){20}}
\put(105,24){\line(4,3){10}}

\put(85,15){\circle*{2}}
\put(100,30){\vector(-1,-1){15}}
\put(100,30){\line(1,1){8}}
\put(100,35){\vector(-3,-4){15}}
\put(85,15){\line(3,4){18}}

\put(166,15){\circle*{2}}
\put(186,15){\vector(-1,0){20}}
\put(186,35){\vector(-1,-1){20}}

\put(166,6){\circle*{2}}
\put(186,6){\vector(-1,0){20}}
\put(186,26){\vector(-1,-1){20}}
\put(186,26){\line(1,1){8}}

\put(0,90){\circle{20}}
\put(0,120){\circle{20}}

\put(30,90){\circle{20}}
\put(30,120){\circle{20}}

\put(80,90){\circle{20}}
\put(80,120){\circle{20}}

\put(110,90){\circle{20}}
\put(110,120){\circle{20}}

\put(160,90){\circle{20}}
\put(160,120){\circle{20}}

\put(190,90){\circle{20}}
\put(190,120){\circle{20}}

\put(-6,86){$A_1$}
\put(-6,118){$A_0$}

\put(74,86){$A_1$}
\put(74,118){$A_0$}

\put(154,86){$A_1$}
\put(154,118){$A_0$}

\put(22,89){\tiny$X_1\backslash A_1$}
\put(22,118){\tiny $X_0\backslash A_0$}

\put(102,87){\tiny$X_1\backslash A_1$}
\put(102,118){\tiny $X_0\backslash A_0$}

\put(182,89){\tiny$X_1\backslash A_1$}
\put(182,118){\tiny $X_0\backslash A_0$}

\put(0,70){$Class\,5'$}
\put(80,70){$Class\,6$}
\put(160,70){$Class\,7'$}
\put(0,-10){$Class\,8'$}
\put(80,-10){$Class\,9'$}
\put(160,-10){$Class\,10'$}
\put(-5,96){\vector(0,1){20}}
\put(2,115){\vector(0,-1){20}}

\put(75,96){\vector(0,1){20}}
\put(82,115){\vector(0,-1){20}}

\put(155,96){\vector(0,1){20}}
\put(162,115){\vector(0,-1){20}}

\put(5,117){\circle*{2}}
\put(5,123){\circle*{2}}
\put(25,117){\vector(-1,0){20}}
\put(25,112){\vector(-4,1){20}}
\put(25,123){\vector(-1,0){20}}
\put(25,128){\vector(-4,-1){20}}

\put(85,112){\circle*{2}}
\put(105,92){\vector(-1,1){20}}
\put(105,97){\vector(-4,3){20}}
\put(105,97){\line(4,-3){8}}

\put(88,116){\circle*{2}}
\put(108,96){\vector(-1,1){20}}
\put(108,101){\vector(-4,3){20}}
\put(108,101){\line(4,-3){8}}

\put(165,116){\circle*{2}}
\put(185,116){\vector(-1,0){20}}
\put(168,122){\circle*{2}}
\put(188,122){\vector(-1,0){20}}

\put(185,96){\vector(-1,1){20}}
\put(188,102){\vector(-1,1){20}}
\put(188,102){\line(1,-1){8}}
\end{picture}
\end{center}
\end{figure}

Where $A$ is a $\lambda-$superatom and $A_i=X_i\cap A$, and\\Class
$5^{'}$ satisfies that $|S_0|-r_0=2$, $|S_1|-r_1=0$, $|T_0|-p=0$ and
$|T_1|-q=0$, and\\ Class $6^{'}$ satisfies that $|S_0|-r_0=0$,
$|S_1|-r_1=0$, $|T_0|-p=0$ and $|T_1|-q=2$, and \\Class $7^{'}$
satisfies that $|S_0|-r_0=1$, $|S_1|-r_1=0$, $|T_0|-p=0$ and
$|T_1|-q=1$,
and\\
Class $8^{'}$ satisfies that $|S_0|-r_0=0$, $|S_1|-r_1=2$,
$|T_0|-p=0$ and $|T_1|-q=0$, and\\Class $9^{'}$ satisfies that
$|S_0|-r_0=0$, $|S_1|-r_1=0$, $|T_0|-p=2$ and $|T_1|-q=0$, and\\
Class $10^{'}$ satisfies that $|S_0|-r_0=0$, $|S_1|-r_1=1$,
$|T_0|-p=1$ and $|T_1|-q=0$, and\\ all of the above digraphs satisfy
that $|A_0|=|A_1|=\delta(X)/2.$

Clearly,\\Class $5^{'}$ is equivalent to Class 5,\\ Class $6^{'}$
is equivalent to Class 9, \\Class $7^{'}$ is equivalent to Class
12  \\Class $8^{'}$ is equivalent to Class 8， \\Class $9^{'}$ is
equivalent to Class 6, and \\ Class $10^{'}$ is equivalent to
Class 13.

\begin{figure}[h,t]
\setlength{\unitlength}{0.46mm}
\begin{center}
\begin{picture}(290,56)

\put(0,10){\circle{20}}
\put(0,40){\circle{20}}

\put(30,10){\circle{20}}
\put(30,40){\circle{20}}

\put(80,10){\circle{20}}
\put(80,40){\circle{20}}

\put(110,10){\circle{20}}
\put(110,40){\circle{20}}

\put(160,10){\circle{20}}
\put(160,40){\circle{20}}

\put(190,10){\circle{20}}
\put(190,40){\circle{20}}

\put(240,10){\circle{20}}
\put(240,40){\circle{20}}

\put(270,10){\circle{20}}
\put(270,40){\circle{20}}
\put(0,-10){$Class\,11'$}
\put(80,-10){$Class\,12'$}
\put(160,-10){$Class\,13'$}
\put(240,-10){$Class\,14'$}

\put(-6,6){$A_1$}
\put(-6,38){$A_0$}

\put(74,6){$A_1$}
\put(74,38){$A_0$}

\put(154,6){$A_1$}
\put(154,38){$A_0$}

\put(234,6){$A_1$}
\put(234,38){$A_0$}

\put(22,9){\tiny$X_1\backslash A_1$}
\put(22,38){\tiny $X_0\backslash A_0$}

\put(102,7){\tiny$X_1\backslash A_1$}
\put(102,38){\tiny $X_0\backslash A_0$}

\put(182,9){\tiny$X_1\backslash A_1$}
\put(182,38){\tiny $X_0\backslash A_0$}

\put(262,9){\tiny$X_1\backslash A_1$}
\put(262,40){\tiny $X_0\backslash A_0$}

\put(-5,16){\vector(0,1){20}}
\put(2,35){\vector(0,-1){20}}

\put(75,16){\vector(0,1){20}}
\put(82,35){\vector(0,-1){20}}

\put(155,16){\vector(0,1){20}}
\put(162,35){\vector(0,-1){20}}

\put(235,16){\vector(0,1){20}}
\put(242,35){\vector(0,-1){20}}
\put(25,36){\vector(-1,0){20}}
\put(25,43){\vector(-1,0){20}}
\put(25,6){\vector(-1,0){20}}
\put(25,13){\vector(-1,0){20}}

\put(105,36){\vector(-1,0){20}}
\put(105,43){\vector(-1,0){20}}
\put(105,35){\vector(-1,-1){22}}
\put(106,32){\vector(-1,-1){22}}

\put(185,6){\vector(-1,0){20}}
\put(185,13){\vector(-1,0){20}}
\put(185,15){\vector(-1,1){20}}
\put(190,16){\vector(-1,1){23}}

\put(268,31){\vector(-1,-1){20}}
\put(268,38){\vector(-1,-1){23}}
\put(265,15){\vector(-1,1){20}}
\put(270,16){\vector(-1,1){23}}
\end{picture}
\end{center}
\end{figure}

Where $A$ is a $\lambda-$superatom and $A_i=X_i\cap A$, and\\ Class
11$^{'}$ satisfies that $|S_0|-r_0=1$, $|S_1|-r_1=1$, $|T_0|-p=0$
and $|T_1|-q=0$, and\\ Class 12$^{'}$ satisfies that $|S_0|-r_0=1$,
$|S_1|-r_1=0$, $|T_0|-p=1$ and $|T_1|-q=0$, and \\Class 13$^{'}$
satisfies that $|S_0|-r_0=0$, $|S_1|-r_1=1$, $|T_0|-p=0$ and
$|T_1|-q=1$,
and\\
Class 14$^{'}$ satisfies that $|S_0|-r_0=0$, $|S_1|-r_1=0$,
$|T_0|-p=1$ and $|T_1|-q=1$, and\\ all of the above digraphs satisfy
that $|A_0|=|A_1|=\delta(X)/2.$
Clearly,\\Class $11^{'}$ is equivalent to Class 11,\\ Class
$12^{'}$ is equivalent to Class 7, \\Class $13^{'}$ is equivalent
to Class 10 , and \\Class $14^{'}$ is equivalent to Class 14.\\\\
All of  the kinds of the special digraphs  of $\mathcal {F}$ are
denoted by $\mathcal
{R}$.\\\\
\textbf{Theorem 4.9} Let $X=MD(G,S_0,S_1,T_0,T_1)$ be max$-\lambda$,
but not super$-\lambda$, if $X$ is neither a directed cycle nor a
cycle and $X$ doesn$^{'}$t belong to $\mathcal {F}$, then $X$
belongs to $\mathcal {L}$ if and only if $X$  belongs to $\mathcal
{R}$\\\\Proof. Necessity. Because $X$ is not super$-\lambda$ and all
the $\lambda-$superatoms contain at least one vertex of $X_0$ and
$X_1$ respectively.
\\\textbf{2.1} $A$ is
a positive $\lambda-$super atom.

Then
$\delta(X)=\lambda(X)=|\omega_X^{+}(A)|=|A_0|(|S_0|-r_0+|T_0|-p)+|A_1|(|S_1|-r_1+|T_1|-q)$.\\Since
$|A|\geq\delta(X)$ and $|A_0|=|A_1|$, we have
$|A_0|=|A_1|\geq\delta(X)/2$. Then \\if
$\lambda(X)=\delta(X)=|\omega_X^{+}(A)|$ only if one of the
following conditions holds.\\\textbf{Case 1} $|S_0|-r_0+|T_0|-p=1$
and $|S_1|-r_1+|T_1|-q=0$.

\texttt{Subcase 1.1}  $|S_0|-r_0=1$, $|T_0|-p=0$, $|S_1|-r_1=0$ and
$|T_1|-q=0$, it is Class 1 or Class $1^{'}$.

\texttt{ Subcase 1.2} $|S_0|-r_0=0$, $|T_0|-p=1$, $|S_1|-r_1=0$ and
$|T_1|-q=0$, it is Class 2 or Class $4^{'}$.\\\textbf{Case 2}
$|S_0|-r_0+|T_0|-p=0$ and $|S_1|-r_1+|T_1|-q=1$.

\texttt{Subcase 2.1} $|S_0|-r_0=0,|T_0|-p=0,|S_1|-r_1=1$ and
$|T_1|-q=0$, it is Class 3 or Class $3^{'}$.

\texttt{Subcase 2.2 } $|S_0|-r_0=0,|T_0|-p=0,|S_1|-r_1=0$ and
$|T_1|-q=1$, it is Class 4 or Class $2^{'}$. \\\textbf{Case 3}
$|S_0|-r_0+|T_0|-p=2$ and $|S_1|-r_1+|T_1|-q=0$.

\texttt{Subcase 3.1} $|S_0|-r_0=2,|T_0|-p=0,|S_1|-r_1=0$ and
$|T_1|-q=0$, it is Class 5 or Class $5^{'}$.

\texttt{Subcase 3.2} $|S_0|-r_0=0,|T_0|-p=2,|S_1|-r_1=0$ and
$|T_1|-q=0$, it is Class 6 or Class $9^{'}$.

\texttt{Subcase 3.3} $|S_0|-r_0=1,|T_0|-p=1,|S_1|-r_1=0$ and
$|T_1|-q=0$, it is Class 7 or Class $12^{'}$. \\\textbf{Case 4}
$|S_0|-r_0+|T_0|-p=0$ and $|S_1|-r_1+|T_1|-q=2$.

\texttt{Subcase 4.1} $|S_0|-r_0=0,|T_0|-p=0,|S_1|-r_1=2$ and
$|T_1|-q=0$, it is Class 8 or Class $8^{'}$.

\texttt{Subcase 4.2} $|S_0|-r_0=0,|T_0|-p=0,|S_1|-r_1=0$ and
$|T_1|-q=2$, it is Class 9 or Class $6^{'}$.

\texttt{Subcase 4.3} $|S_0|-r_0=0,|T_0|-p=0,|S_1|-r_1=1$ and
$|T_1|-q=1$, it is Class 10 or Class $13^{'}$. \\\textbf{Case 5}
$|S_0|-r_0+|T_0|-p=1$ and $|S_1|-r_1+|T_1|-q=1$.

Similarly, we can deduce that under this case, it is Class 11,
Class 12, Class 13, Class 14, Class $7^{'}$, Class $10^{'}$,
Class $11^{'}$ or Class $14^{'}$.\\Sufficiency. Clearly. $\Box$\\
\\\textbf{Proposition 4.10} $X=MD(G,S_0,S_1,T_0,T_1)$ is a strongly connected mixed
Cayley digraph, and $X$ belongs to $\mathcal {L}$, then $X$ is
Class 1 or Class $1^{'}$ if and only if \\(1)There exists a
non-trivial proper subgroup $H$ of $G$ and $S_0$ contains an element
$s_0$ such that

$<S_0\cup \{1_G\}\backslash\{s_0\}>\leq H$ and $|H|=\delta(X)$,
and\\(2) There is an element $t_0\in T_0$ such that

$G_1=<S_1>\leq t_0Ht_0^{-1}$, $T_1^{-1}t_0\subseteq H$ and
$t_0^{-1}T_0\subseteq H$.\\\\Proof. Necessity. Because the class of
Class 1 is equivalent to the class of Class $1^{'}$, without loss
of generality, we set $X$ belongs to the class of Class 1, then
Assume $(1_{G},0)\in A_0$, by lemma 4.5, $H_0\leq G$. Let $H=H_0$,
then under this situation we can achieve the following results
easily,\\(i)
$\lambda(X)=|\omega_X^{+}(A)|=|A_0|=|H_0|=|H|=\delta(X)$, since
$|S_0|-r_0=1$, $|T_0|-p=0$, $|S_1|-r_1=0$, and $|T_1|-q=0$,\\(ii)
$<S_0\cup\{1_G\}\setminus\{s_0\}>\leq H_0=H$, since $|S_0|-r_0=1$.\\
By proposition 4.7, $H_1=t_0H_0$ for some $t_0\in T_0$ and
$X_1=\cup_{i=1}^{k}(t_0H_0g_i)\times\{1\}$,\\ where $t_0H_0g_i\cap
t_0H_0g_j\neq \varnothing$ if and only if $i=j$ for $1\leq i,j\leq
k$. \\Assume that $(1_G,1)\in (t_0H_0g_s)\times\{1\}$, then we can
deduce that

$t_0H_0g_s\leq G$ and $g_s=h_0^{-1}t_0^{-1}$, where $h_0\in
H_0$.\\Since $|S_1|-r_1=0$, we get $G_1=<S_1>\leq
t_0H_0g_s=t_0H_0h_0^{-1}t_0^{-1}=t_0H_0t_0^{-1}=t_0Ht_0^{-1}$.\\Since
$|T_0|-p=0$ and $|T_1|-q=0$, then

$T_0H_0\subseteq H_1$ and $T_1^{-1}H_1\subseteq H_0$, \\so
$T_0H_0\subseteq t_0H_0$ and $T_{1}^{-1}t_0H_0\subseteq H_0$.\\ It
means that $t_0^{-1}T_0\subseteq H_0=H$  and $T_1^{-1}t_0\subseteq
H_0=H$  for some $t_0\in T_0$.\\Sufficiency, set $A=H\times
\{0\}\cup (t_0H)\times \{1\}$,
\\because $<S_0\cup \{1_G\}\setminus \{s_0\}>\leq H$, we can
get $|S_0|-r_0=1$. \\Similarly, because

$G_1=<S_1>\leq t_0Ht_0^{-1}$, $T_1^{-1}t_0\subseteq H$ and
$t_0^{-1}T_0\subseteq H$,
\\we can get that

$|S_1|-r_1=0$, $|T_0|-p=0$ and  $|T_1|-q=0$,\\ so
$\lambda(X)=|\omega^{+}(A)|=|H|=\delta(X)$, and $A$ is not
nontrivial. $\Box$\\\\Analogously, we can get the following
proposition from 4.11 to 4.23.
\\\\\textbf{Proposition 4.11} $X=MD(G,S_0,S_1,T_0,T_1)$ is a strongly connected mixed
Cayley digraph, and $X$ belongs to $\mathcal {L}$, then $X$ is
Class 2 or Class $4^{'}$ if and only if \\$(1)$ There exists a
non-trivial proper subgroup $H$ of $G$ such that

$G_0=<S_0>\leq H$ and $|H|=\delta(X)$, and\\$(2)$ There are two
distinct elements $t_0,t_0^{'}\in T_0$ such that

$G_1=<S_1>\leq t_0Ht_0^{-1},T_1^{-1}t_0\subseteq H,$

$t_0^{'}H\cap t_0H=\varnothing$ and
$t_0^{-1}(T_0\setminus\{t_0^{'}\})\subseteq H$.
$\Box$\\\\\textbf{Proposition 4.12} $X=MD(G,S_0,S_1,T_0,T_1)$ is a
strongly connected mixed Cayley digraph, and $X$ belongs to
$\mathcal {L}$, then $X$ is Class 3 or Class $3^{'}$ if and only
if \\(1) There is a non-trivial proper subgroup $H$ of $G$ and some
element $s_1\in S_1$ such that

$<S_1\cup\{1_G\}\backslash \{s_1\}>\leq H$  and $|H|=\delta(X)$,
and\\(2) There is an element $t_1\in T_1$ such that

$G_0=<S_0>\leq t_1^{-1}Ht_1,T_0t_1^{-1}\subseteq H$ and
$t_1T_1^{-1}\subseteq H$. $\Box$\\\\\textbf{Proposition 4.13}
$X=MD(G,S_0,S_1,T_0,T_1)$ is a strongly connected mixed Cayley
digraph, and $X$ belongs to $\mathcal {L}$, then $X$ is Class 4 or
Class $2^{'}$ if and only if \\(1) There exists a non-trivial
proper subgroup $H$ of $G$ such that

$G_1=<S_1>\leq H$ and $|H|=\delta(X)$, and \\(2) There are two
distinct elements $t_1,t_1^{'}\in T_1$ such that

$G_0=<S_0>\leq t_1^{-1}Ht_1,T_0t_1^{-1}\subseteq H$,

${t_1^{'}}^{-1}H\cap t_1^{-1}H=\varnothing$ and
$t_1(T_1^{-1}\setminus\{{t_1^{'}}^{-1}\})\subseteq H$.
$\Box$\\\\\textbf{Proposition 4.14} $X=MD(G,S_0,S_1,T_0,T_1)$ is a
strongly connected mixed Cayley digraph, and $X$ belongs to
$\mathcal {L}$, then $X$ is Class 5 or Class $5^{'}$ if and only
if
\\(1) There exists a non-trivial proper subgroup $H$ of $G$ and
$S_0$ contains two distinct elements $s_0,s_0^{'}$ such that

$<S_0\cup\{1_G\}\backslash \{s_0,s_0^{'}\}>\leq H$ and
$|H|=\delta(X)/2$, and\\(2) There is an element $t_0\in T_0$ such
that

$G_1=<S_1>\leq t_0Ht_0^{-1}$, $T_1^{-1}t_0\subseteq H$ and
$t_0^{-1}T_0\subseteq H$. $\Box$\\\\\textbf{Proposition 4.15}
$X=MD(G,S_0,S_1,T_0,T_1)$ is a strongly connected mixed Cayley
digraph, and $X$ belongs to $\mathcal {L}$, then $X$ is Class 6 or
Class $9^{'}$ if and only if \\(1) There exists a non-trivial
subgroup $H$ of $G$ such that

$G_0=<S_0>\leq H$ and $|H|=\delta(X)/2$, and\\(2) There are three
distinct elements $t_0,t_0^{'},t_0^{''}\in T_0$ such that

$G_1=<S_1>\leq t_0Ht_0^{-1}$, $T_1^{-1}t_0\subseteq H$,
$t_0^{'}H\cap t_0H=\varnothing$,

$t_0^{''}H\cap t_0H=\varnothing$ and  $t_0^{-1}(T_0\setminus
\{t_0^{'}$, $t_0^{''}\})\subseteq H$. $\Box$\\\\\textbf{Proposition
4.16} $X=MD(G,S_0,S_1,T_0,T_1)$ is a strongly connected mixed Cayley
digraph, and $X$ belongs to $\mathcal {L}$, then $X$ is Class 7 or
Class $12^{'}$ if and only if\\(1) There exists a non-trivial
proper subgroup $H$ of $G$,and $S_0$ contains an element $s_0$ such
that

$<S_0\cup\{1_G\}\backslash\{s_0\}>\leq H$ and $|H|=\delta(X)/2$,
and\\(2) There are two distinct elements $t_0,t_0^{'}\in T_0$ such
that

$G_1=<S_1>\leq t_0Ht_0^{-1}$, $T_1^{-1}t_0\subseteq H$,

$t_0^{'}H\cap t_0H=\varnothing$ and
$t_0^{-1}(T_0\setminus\{t_0^{'}\})\subseteq H$. $\Box$
\\\\\textbf{Proposition 4.17} $X=MD(G,S_0,S_1,T_0,T_1)$ is a
strongly connected \emph{mixed Cayley digraph }, and $X$ belongs to
$\mathcal {L}$, then $X$ is Class 8 or Class $8^{'}$ if and only
if\\(1) There is a non-trivial subgroup $H$ of $G$ and  some
$s_1,s_1^{'}\in S_1$ such that

$<S_1\cup \{1_G\}\backslash\{s_1,s_1^{'}\}>\leq H$  and
$|H|=\delta(X)/2$, and \\(2) There is an element $t_1\in T_1$ such
that

$G_0=<S_0>\leq t_1^{-1}Ht_1$, $T_0t_1^{-1}\subseteq H$ and
$t_1T_1^{-1}\subseteq H$. $\Box$\\\\\textbf{Proposition 4.18}
$X=MD(G,S_0,S_1,T_0,T_1)$ is a strongly connected mixed Cayley
digraph, and $X$ belongs to $\mathcal {L}$, then $X$ is Class 9 or
Class $6^{'}$ if and only if\\(1) There is a non-trivial proper
subgroup $H$ of $G$ such that

$G_1=<S_1>\leq H$ and $|H|=\delta(X)/2$, and\\(2) There are there
distinct elements $t_1$, $t_1^{'}$, $t_1^{''}\in T_1$ such that

$G_0=<S_0>\leq t_1^{-1}Ht_1,T_0t_1^{-1}\subseteq H$, $t_1^{-1}H\cap
{t_1^{'}}^{-1}H=\varnothing$,

$t_1^{-1}H\cap {t_1^{''}}^{-1}H=\varnothing$ and
$t_1(T_1^{-1}\setminus\{{t_1^{'}}^{-1}$,
${t_1^{''}}^{-1}\})\subseteq H$. $\Box$\\\\\textbf{Proposition 4.19}
$X=MD(G,S_0,S_1,T_0,T_1)$ is a strongly connected mixed Cayley
digraph, and $X$ belongs to $\mathcal {L}$, then $X$ is Class 10 or
Class $13^{'}$ if and only if\\(1) There is a non-trivial proper
subgroup $H$ of $G$ and some element $s_1\in S_1$ such that

$<S_1\cup\{1_G\}\backslash \{s_1\}>\leq H$  and $|H|=\delta(X)/2$,
and\\(2) There are two distinct elements $t_1$, $t_1^{'}\in T_1$
such that

$G_0=<S_0>\leq t_1^{-1}Ht_1$, $T_0t_1^{-1}\subseteq H$,

${t_1^{'}}^{-1}H\cap t_1^{-1}H=\varnothing$ and
$t_1(T_1^{-1}\setminus\{{t_1^{'}}^{-1}\})\subseteq H$.
$\Box$\\\\\textbf{Proposition 4.20} $X=MD(G,S_0,S_1,T_0,T_1)$ is a
strongly connected mixed Cayley digraph, and $X$ belongs to
$\mathcal {L}$, then $X$ is Class 11 or Class $11^{'}$ if and only
if\\(1) There is a non-trivial proper subgroup $H$ of $G$ and $S_0$
contains an element $s_0$ such that

$<S_0\cup\{1_G\}\backslash \{s_0\}>\leq H$ and $|H|=\delta(X)/2$,
and\\(2) There is an element $t_0\in T_0$ and an element $s_1\in
S_1$ such that

$<S_1\cup\{1_G\}\backslash \{s_1\}>\leq t_0Ht_0^{-1}$ and
$T_1^{-1}t_0$, $t_0^{-1}T_0\subseteq H$.
$\Box$\\\\\textbf{Proposition 4.21} $X=MD(G,S_0,S_1,T_0,T_1)$ is a
strongly connected mixed Cayley digraph, and $X$ belongs to
$\mathcal {L}$, then $X$ is Class 12 or Class $7^{'}$ if and only
if
\\(1) There is a non-trivial proper subgroup $H$ of $G$ and $S_0$
contains an element $s_0$ such that

$<S_0\cup\{1_G\}\backslash \{s_0\}>\leq H$ and $|H|=\delta(X)/2$,
and \\(2) There is an element $t_0\in T_0$ and an element $t_1\in
T_1$ such that

$G_1=<S_1>\leq t_0Ht_0^{-1}$, $t_0^{-1}T_0\subseteq H$ ,

$t_1^{-1}t_0\notin H$ and
$(T_1^{-1}\setminus\{t_1^{-1}\})t_0\subseteq H$. $\Box$
\\\\\textbf{Proposition 4.22} $X=MD(G,S_0,S_1,T_0,T_1)$ is a
strongly connected mixed Cayley digraph, and $X$ belongs to
$\mathcal {L}$, then $X$ is Class 13 or Class $10^{'}$ if and only
if  \\(1) There is a non-trivial proper subgroup $H$ of $G$ and some
element $s_1\in S_1$ such that

$<S_1\cup\{1_G\}\backslash\{s_1\}>\leq H$  and $|H|=\delta(X)/2$,
and\\(2) There is an element $t_1\in T_1$ and an element $t_0\in
T_0$ such that

$G_0=<S_0>\leq t_1^{-1}Ht_1$, $t_1T^{-1}\subseteq H$,

$t_0t_1^{-1}\notin H$ and $(T_0\setminus\{t_0\})t_1^{-1}\subseteq
H$. $\Box$\\\\\textbf{Proposition 4.23} $X=MD(G,S_0,S_1,T_0,T_1)$ is
a strongly connected  mixed Cayley digraph, and $X$ belongs to
$\mathcal {L}$, then $X$ is Class 14 or Class $14^{'}$ if and only
if  \\(1) There is an non-trivial proper subgroup $H$ of $G$ such
that

$G_0=<S_0>\leq H$ and $|H|=\delta(X)/2$, and\\(2) There are there
distinct elements $t_0,t_0^{'}\in T_0,t_1\in T_1$ such that

$G_1=<S_1>\leq t_0Ht_0^{-1}$, $t_1^{-1}t_0\notin H$,
$(T_1^{-1}\setminus \{t_1^{-1}\})t_0\subseteq H$,

$t_0^{'}H\cap t_0^{'}H=\varnothing$ and
$t_0^{-1}(T_0\setminus\{t_0^{'}\})\subseteq H$. $\Box$\\\\From the
above discussion, we get the following theorem.\\\\\textbf{Theorem
4.24} Let $X=MD(G,S_0,S_1,T_0,T_1)$ be max$-\lambda$, if $X$ is
neither a directed cycle nor a cycle and $X$ doesn$^{'}$t belong to
$F$, then $X$ is not super$-\lambda$ if and only if $X$ satisfies
one of the following conditions:\\(1) $|T_0|=1$ or
$|T_1|=1$,$1\leq|S_0|\leq|S_1|$ and $S_0\cup\{1_G\}\leq G$.\\(2)
$|T_0|=1$ or $|T_1|=1$,$1\leq|S_1|\leq|S_0|$ and $S_1\cup\{1_G\}\leq
G$.\\(3) (3.1) There exists a non-trivial proper subgroup $H$ of $G$
and $S_0$ contains an element $s_0$ such that

$<S_0\cup \{1_G\}\backslash\{s_0\}>\leq H$ and $|H|=\delta(X)$, and

(3.2) There is an element $t_0\in T_0$ such that

$G_1=<S_1>\leq t_0Ht_0^{-1}$, $T_1^{-1}t_0\subseteq H$ and
$t_0^{-1}T_0\subseteq H$.\\(4) $(4.1)$ There exists a non-trivial
proper subgroup $H$ of $G$ such that

$G_0=<S_0>\leq H$ and $|H|=\delta(X)$, and

$(4.2)$ There are two distinct elements $t_0,t_0^{'}\in T_0$ such
that

$G_1=<S_1>\leq t_0Ht_0^{-1},T_1^{-1}t_0\subseteq H,$

$t_0^{'}H\cap t_0H=\varnothing$ and
$t_0^{-1}(T_0\setminus\{t_0^{'}\})\subseteq H$.\\(5) (5.1) There is
a non-trivial proper subgroup $H$ of $G$ and  some element $s_1\in
S_1$ such that

$<S_1\cup\{1_G\}\backslash \{s_1\}>\leq H$  and $|H|=\delta(X)$, and

(5.2) There is an element $t_1\in T_1$ such that

$G_0=<S_0>\leq t_1^{-1}Ht_1,T_0t_1^{-1}\subseteq H$ and
$t_1T_1^{-1}\subseteq H$. \\(6) (6.1) There exists a non-trivial
proper subgroup $H$ of $G$ such that

$G_1=<S_1>\leq H$ and $|H|=\delta(X)$, and

(6.2) There are two distinct elements $t_1,t_1^{'}\in T_1$ such that

$G_0=<S_0>\leq t_1^{-1}Ht_1,T_0t_1^{-1}\subseteq H$,

${t_1^{'}}^{-1}H\cap t_1^{-1}H=\varnothing$ and
$t_1(T_1^{-1}\setminus\{{t_1^{'}}^{-1}\})\subseteq H$. \\(7) (7.1)
There exists a non-trivial proper subgroup $H$ of $G$ and $S_0$
contains two distinct

elements $s_0,s_0^{'}$ such that

$<S_0\cup\{1_G\}\backslash \{s_0,s_0^{'}\}>\leq H$ and
$|H|=\delta(X)/2$, and

(7.2) There is an element $t_0\in T_0$ such that

$G_1=<S_1>\leq t_0Ht_0^{-1}$,$T_1^{-1}t_0\subseteq H$ and
$t_0^{-1}T_0\subseteq H$.  \\(8) (8.1) There exists a non-trivial
subgroup $H$ of $G$ such that

$G_0=<S_0>\leq H$ and $|H|=\delta(X)/2$, and

(8.2) There are three distinct elements $t_0$, $t_0^{'}$,
$t_0^{''}\in T_0$ such that

$G_1=<S_1>\leq t_0Ht_0^{-1}$, $T_1^{-1}t_0\subseteq H$,
$t_0^{'}H\cap t_0H=\varnothing$,

$t_0^{''}H\cap t_0H=\varnothing$ and $t_0^{-1}(T_0\setminus
\{t_0^{'},t_0^{''}\})\subseteq H$.
\\(9) (9.1) There exists a non-trivial proper subgroup $H$ of
$G$, and $S_0$ contains an element

$s_0$ such that

$<S_0\cup\{1_G\}\backslash\{s_0\}>\leq H$ and $|H|=\delta(X)/2$,
and

(9.2) There are two distinct elements $t_0,t_0^{'}\in T_0$ such that

$G_1=<S_1>\leq t_0Ht_0^{-1}$, $T_1^{-1}t_0\subseteq H$,

$t_0^{'}H\cap t_0H=\varnothing$ and
$t_0^{-1}(T_0\setminus\{t_0^{'}\})\subseteq H$. \\(10) (10.1) There
is a non-trivial subgroup $H$ of $G$ and some $s_1,s_1^{'}\in S_1$
such that

$<S_1\cup \{1_G\}\backslash\{s_1,s_1^{'}\}>\leq H$ and
$|H|=\delta(X)/2$, and

(2) There is an element $t_1\in T_1$ such that

$G_0=<S_0>\leq t_1^{-1}Ht_1$, $T_0t_1^{-1}\subseteq H$ and
$t_1T_1^{-1}\subseteq H$. \\(11) (11.1) There is a non-trivial
proper subgroup $H$ of $G$ such that

$G_1=<S_1>\leq H$ and $|H|=\delta(X)/2$, and

(11.2) There are there distinct elements $t_1$, $t_1^{'}$,
$t_1^{''}\in T_1$ such that

$G_0=<S_0>\leq t_1^{-1}Ht_1$, $T_0t_1^{-1}\subseteq H$,
$t_1^{-1}H\cap {t_1^{'}}^{-1}H=\varnothing$,

$t_1^{-1}H\cap {t_1^{''}}^{-1}H=\varnothing$ and
$t_1(T_1^{-1}\setminus\{{t_1^{'}}^{-1}$,
${t_1^{''}}^{-1}\})\subseteq H$. \\(12)(12.1) There is a non-trivial
proper subgroup $H$ of $G$ and some element $s_1\in S_1$ such that

$<S_1\cup\{1_G\}\backslash \{s_1\}>\leq H$ and $|H|=\delta(X)/2$,
and

(12.2) There are two distinct elements $t_1$, $t_1^{'}\in T_1$ such
that

$G_0=<S_0>\leq t_1^{-1}Ht_1$, $T_0t_1^{-1}\subseteq H$,

${t_1^{'}}^{-1}H\cap t_1^{-1}H=\varnothing$ and
$t_1(T_1^{-1}\setminus\{{t_1^{'}}^{-1}\})\subseteq H$. \\(13) (13.1)
There is a non-trivial proper subgroup $H$ of $G$ and $S_0$ contains
an element

 $s_0$ such that

$<S_0\cup\{1_G\}\backslash \{s_0\}>\leq H$ and $|H|=\delta(X)/2$,
and

(13.2) There is an element $t_0\in T_0$ and an element $s_1\in S_1$
such that

$<S_1\cup\{1_G\}\backslash \{s_1\}>\leq t_0Ht_0^{-1}$, $T_1^{-1}t_0$
and $t_0^{-1}T_0\subseteq H$. \\(14) (14.1) There is a non-trivial
proper subgroup $H$ of $G$ and $S_0$ contains an element

$s_0$ such that

$<S_0\cup\{1_G\}\backslash \{s_0\}>\leq H$ and $|H|=\delta(X)/2$,
and

(14.2) There is an element $t_0\in T_0$ and an element $t_1\in T_1$
such that

$G_1=<S_1>\leq t_0Ht_0^{-1}$, $t_0^{-1}T_0\subseteq H$,

$t_1^{-1}t_0\notin H$ and
$(T_1^{-1}\setminus\{t_1^{-1}\})t_0\subseteq H$. \\(15) (15.1) There
is a non-trivial proper subgroup $H$ of $G$ and some element $s_1\in
S_1$ such that

$<S_1\cup\{1_G\}\backslash\{s_1\}>\leq H$ and $|H|=\delta(X)/2$, and

(15.2) There is an element $t_1\in T_1$ and an element $t_0\in T_0$
such that

$G_0=<S_0>\leq t_1^{-1}Ht_1$, $t_1T^{-1}\subseteq H$,

$t_0t_1^{-1}\notin H$ and $(T_0\setminus\{t_0\})t_1^{-1}\subseteq
H$. \\(16) (16.1) There is an non-trivial proper subgroup $H$ of $G$
such that

$G_0=<S_0>\leq H$ and $|H|=\delta(X)/2$, and

(16.2) There are there distinct elements $t_0$, $t_0^{'}\in T_0$,
$t_1\in T_1$ such that

$G_1=<S_1>\leq t_0Ht_0^{-1},t_1^{-1}t_0\notin H$,
$(T_1^{-1}\setminus \{t_1^{-1}\})t_0\subseteq H$,

$t_0^{'}H\cap
t_0^{'}H=\varnothing$ and $t_0^{-1}(T_0\setminus\{t_0^{'}\})\subseteq H$. $\Box$\\

So we can conclude that the strongly connected mixed Cayley digraph
is max$-\lambda$ and super$-\lambda$ but a few exceptions.


\newpage
\begin{thebibliography}{99}
\bibitem{Bondy} J.A. Bondy, U.S.R. Murty. Graph Theory with Applications,
 North-Holland, New York, 1976.
\bibitem{Chen} J.Y.Chen, Jixiang Meng, Super edge-connectivity of
mixed Cayley graph, Science Direct,2007.
\bibitem{Lu} Z.P.Lu, On the Automorphism Groups of Bi-Cayley Graphs,
Acta.Sci.Natu. Universitatics Pekinensis, 39(2003)
\bibitem{Mader} M.Mader, Minimale n-fach Kantenzusammenhangenden
Graphen, Math. Ann. 191(1971)21-28
\bibitem{Meng} J.X.Meng, Connectivity of vertex and edge transitive graphs, Discrete
Appl.Math.127(2003)601-603.
\bibitem{Meng} J.X.Meng, Optimally super-edge-connected
transitive graphs,Discrete Math.260(2003)239-248.
\bibitem{MJX} J.X.Meng, Sper-connectivity of Vertex-Transitive Bi-Cayley Graphs and Bipartite Cayley Graph.
\bibitem{Tindell} R. Tindell, Connectivity of Cayley
digraphs, in:D.Z.Du, D.F.Hsu(Eds.), Combinatorial Network
Theory,Klumer,Dordrecht,1996,pp.41-46.
\bibitem{Watkins} M.E.Watkins, Connectivity of transitive graphs,
J.Comb. Theory 8*1970)23-29.
\bibitem{Xu} M.Y.Xu, Introduction of Finite Group, Vol $\amalg$,Science
Press,Beijing,1999,
 384-386.
\end{thebibliography}
\end{document}